%% file: text17.tex
\documentclass[11pt]{nyj}
\usepackage{graphicx,color}
\usepackage{fourier}
\graphicspath{{pictures/}}

\swapnumbers
\theoremstyle{plain}
\newtheorem{theorem}{Theorem}[section]
\newtheorem{proposition}[theorem]{Proposition}
\newtheorem{lemma}[theorem]{Lemma}
\newtheorem{corollary}[theorem]{Corollary}

\theoremstyle{definition}
\newtheorem{definition}[theorem]{Definition}
\newtheorem{question}[theorem]{Question}

\theoremstyle{remark}
\newtheorem{remark}[theorem]{Remark}

\newcommand{\hyp}{\nobreakdash-\hspace{0pt}}
\def\BZ{\mathbb Z}

\def\Mtau{M^{\tau}}
\def\Ktau{K^{\tau}}
\def\Ttau{T^{\tau}}

\makeatletter
  \let\c@theorem=\c@subsection
  \let\c@figure=\c@subsection
  \let\p@figure=\p@subsection
  \let\thefigure=\thesubsection
  \let\cl@figure=\cl@subsection
  \let\c@table=\c@subsection
  \let\p@table=\p@subsection
  \let\thetable=\thesubsection
  \let\cl@table=\cl@subsection
  \let\c@equation=\c@subsection
  \let\p@equation=\p@subsection
  \let\theequation=\thesubsection
  \let\cl@equation=\cl@subsection
\makeatother

\newcommand{\printdata}[5]{%
  \begin{center}
    \rotatebox{90}{\mbox{\begin{minipage}[c][11.5cm][c]{18.5cm}
          \begin{center}
            {\large Almost mutant pair: $\mathbf{#3}$ and $\mathbf{#4}$}\\
          \end{center}{
            \vspace{0.5cm}
            \renewcommand{\arraystretch}{1.2}
            \begin{tabular}{lll}
              Volume: \csname HVolume_#1\endcsname & Signature: $#5$  & Alexander: \csname Alexander_#1\endcsname \\ 
          \multicolumn{3}{l}{Jones: \csname Jones_#1\endcsname} \\ 
        \end{tabular}\\
        \renewcommand{\arraystretch}{1}
        \\
          \raisebox{1.5cm}{\qquad HOMFLY-PT:} \hspace{-2cm} \csname HOMFLYPT_#1\endcsname \qquad  \qquad\raisebox{1.5cm}{Kauffman:} {\footnotesize \csname Kauffman_#1\endcsname} \\
          \renewcommand{\arraystretch}{1.5}
          \begin{tabular}{rl}
            \multicolumn{2}{l}{Khovanov Homology for $\mathbf{#3}$:}\\ 
            ranks: & \csname KhovanovQ_#1\endcsname \\ 
            2-torsion: & \csname KhovanovT_#1\endcsname \\
            \multicolumn{2}{l}{Khovanov Homology for $\mathbf{#4}$:}\\
          ranks:&\csname KhovanovQ_#2\endcsname\\
          2-torsion:&\csname KhovanovT_#2\endcsname\\
        \end{tabular}
        \renewcommand{\arraystretch}{1}
      }\end{minipage}
  }}
\end{center}
\vfill 
\pagebreak
}

\title[Behavior of knot invariants under  genus $2$ mutation]{
Behavior of knot invariants under \\ genus $2$ mutation}

\author[Dunfield]{Nathan M.~Dunfield}
\address{Dept. of Mathematics, MC-382\\
  University of Illinois \\ 
  Urbana, IL 61801,  USA}
\email{nathan@dunfield.info}
\urladdr{http://dunfield.info}

\author[Garoufalidis]{Stavros Garoufalidis}
\address{School of Mathematics \\
          Georgia Institute of Technology \\
          Atlanta, GA 30332-0160, USA}
\email{stavros@math.gatech.edu}
\urladdr{http://www.math.gatech.edu/$\sim$stavros }

\author[Shumakovitch]{Alexander Shumakovitch}
\address{George Washington University \\
          Department of Mathematics \\
          1922 F Street, NW \\
          Washington, DC 20052, USA}
\email{shurik@gwu.edu}

\author[Thistlethwaite]{Morwen Thistlethwaite}
\address{Department of Mathematics \\
          The University of Tennessee \\
          Knoxville, TN 37996-1300, USA}
\email{morwen@math.utk.edu}
\urladdr{http://www.math.utk.edu/$\sim$morwen}

\thanks{N.D.~was partially supported by the supported by the Sloan Foundation.  N.D.~and S.G.~were partially supported by the U.S.~N.S.F. \\
This article appeared as \emph{New York J. Math.}, \textbf{16} (2010) 99-123.} 


\begin{document}

\begin{abstract}
  Genus 2 mutation is the process of cutting a 3-manifold along an
  embedded closed genus 2 surface, twisting by the hyper-elliptic
  involution, and gluing back.  This paper compares genus 2 mutation
  with the better-known Conway mutation in the context of knots in the
  3-sphere.  Despite the fact that any Conway mutation can be achieved
  by a sequence of at most two genus 2 mutations, the invariants that are
  preserved by genus 2 mutation are a proper subset of those preserved
  by Conway mutation.  In particular, while the Alexander and Jones
  polynomials are preserved by genus 2 mutation, the HOMFLY-PT
  polynomial is not.  In the case of the $sl_2$-Khovanov homology,
  which may or may not be invariant under Conway mutation, we give an
  example where genus 2 mutation changes this homology.  Finally,
  using these techniques, we exhibit examples of knots with the same
  same colored Jones polynomials, HOMFLY-PT polynomial, Kauffman
  polynomial, signature and volume, but different Khovanov homology.
\end{abstract}

\maketitle

\tableofcontents

\section{Introduction}
\label{sec.intro}

In the 1980s, a plethora of new knot invariants were discovered,
following the discovery of the Jones polynomial \cite{J}. These
powerful invariants were by construction {\em chiral}, i.e.~they were
often able to distinguish knots from their mirrors, as opposed to many
of their classical counterparts.  Soon after the appearance of these
new quantum invariants of knots, many people studied their behavior
under other kinds of involutions, and in particular under mutation.
Chmutov, Duzhin, Lando, Lickorish, Lipson, Morton, Traczyk and others
pioneered the behavior of the quantum knot invariants under mutation;
see \cite{CDL,LL,MC,MR1} and references therein. The quantum invariants
come in two flavors: rationally valued Vassiliev invariants, and
polynomially valued exact invariants (such as the Jones, HOMFLY,
Kauffman, Alexander polynomials), see \cite{Tu2}.  Later on, abelian
group valued invariants were constructed by Khovanov \cite{Kh}.

\begin{figure}[t]
  \centerline{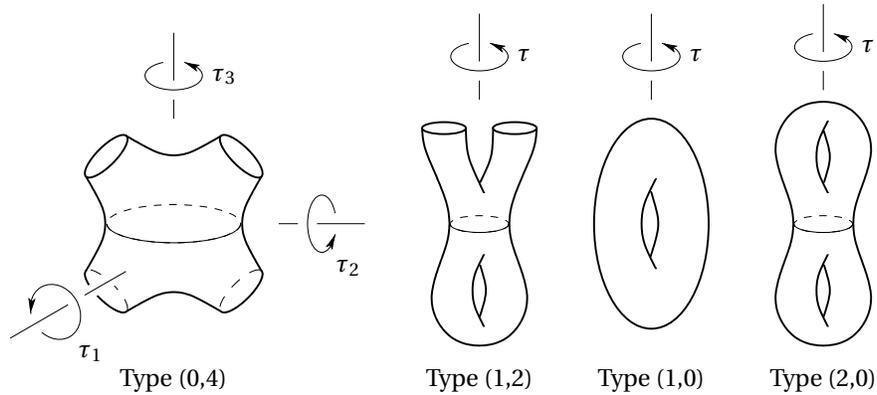}
  \caption{Symmetric surfaces of types $(0,4)$, $(1,2)$, $(1,0)$, and 
    $(2,0)$ and their involutions.  There are also symmetric surfaces of
    type $(1,1)$ and $(0,3)$ that are not pictured, since we will not
    need them here.}\label{fig:symmetricsurface}
\end{figure}

Here, we study the behavior of classical and quantum invariants of
knots in $S^3$ under mutation, building on the above mentioned work.
The notion of mutation was introduced by Conway in~\cite{Co}, and has
been used extensively in various generalized forms.  Let us start by
explaining what we mean by mutation.  Roughly, mutation is modifying a
3-manifold by cutting it open along a certain kind of embedded
surface, and then regluing in a different way.  More precisely,
consider one of the surfaces $F$ from
Figure~\ref{fig:symmetricsurface}, together with the specified
involution $\tau$; we will call the pair $(F,\tau)$ a {\em symmetric
  surface}.
Suppose $F$ is a symmetric surface properly embedded in a compact
orientable 3-manifold. The \emph{mutant} of $M$ along $F$ is the
result of cutting $M$ open along $F$, and then regluing the two copies
of $F$ by the involution $\tau$.  The mutant manifold is denoted
$M^\tau$, and the operation is called \emph{mutation}.  When we want to
distinguish the topological type of $F$, we refer to $(g,s)$-mutation
where $g$ is the genus and $s$ is the number of boundary components.

The involutions used in mutation have very special properties, e.g.~
if $\gamma$ is a non--boundary-parallel simple closed curve, then 
$\tau(\gamma)$
is isotopic to $\gamma$ (neglecting orientations).  As a result, while
mutation is typically violent enough to change the global topology of
$M$, it is simultaneously subtle enough that many invariants do not
change.  Studying this phenomenon has enriched our understanding of a
number of invariants, be they classical, quantum, or geometric.

When studying knots in $S^3$, the most natural type of mutation is
$(0,4)$\hyp mutation, which has a simple interpretation in terms of a knot
diagram, and is known to preserve a wide range of invariants.  Here,
we study the effects of $(2,0)$-mutation on knots in $S^3$.  By this,
we mean the following.  If $F$ is a closed 2-surface in $S^3$, then
the mutant $(S^3)^\tau$ is always homeomorphic to $S^3$ (see
Section~\ref{subsec.knotmut}).  Thus if $K$ is a knot in $S^3$ which
is disjoint from $F$, it makes sense to talk about its mutant
$K^\tau$.

In this context, $(2,0)$-mutation is the most general type: any of the
above mutations can be achieved by a sequence of at most two
$(2,0)$-mutations (see Lemma~\ref{lem.20suffice} below).  Given that
any $(0,4)$-mutation can be implemented in this way, you might expect
that an invariant unchanged by $(0,4)$-mutation would also be
preserved by $(2,0)$-mutations.  It turns out that this is not the
case, as you can see from the following table; with the possible
exception of Khovanov homology, all of the invariants listed there are
preserved by $(0,4)$-mutation.
\begin{table}[htb]
  \begin{tabular}{c|c} 
    Preserved by $(2,0)$-mutation  &  Changed by $(2,0)$-mutation \\[5pt] \hline \\[-5pt]

        Hyperbolic volume/Gromov norm \\ of the knot exterior  & HOMFLY-PT polynomial \\ 
        Alexander polynomial and \\ generalized signature &  $sl_2$-Khovanov Homology \\
    Colored Jones polynomials   & \\ 
  \end{tabular}
  
  \vspace{0.5cm}

  \caption{Summary of known results on genus 2 mutation.}
  \label{table.main}
\end{table}

The results on the left-hand side are due either entirely or in large
part to Ruberman~\cite{Ru}, Cooper-Lickorish~\cite{CL} and
Morton-Traczyk~\cite{MT}, see below for details; the results on the
right are new.  One way of interpreting these results might be that
the invariants on the left are more tied to the topology of $S^3 \setminus
K$, whereas those on the right are more ``diagrammatic'' and tied to
combinatorics of knot projections.  (Of course, this must be taken
with more than a grain of salt, since knots are determined by their
complements~\cite{GordonLuecke}.)  The presence of the colored Jones
polynomials among the more ``topological'' invariants is not so
surprising given their connections to purely geometric/topological
invariants in the context of the Volume Conjecture (e.g.~the results
of \cite{GT}).  Indeed, one of our original motivations for this work
was to better understand the Volume Conjecture, which proposes a
relationship between the colored Jones polynomials and the hyperbolic
volume.  The fact that both the colored Jones polynomials and
hyperbolic volume are preserved by $(2,0)$-mutation is positive
evidence for this conjecture.

\begin{figure}[htb]
\centerline{\includegraphics[scale=.3]{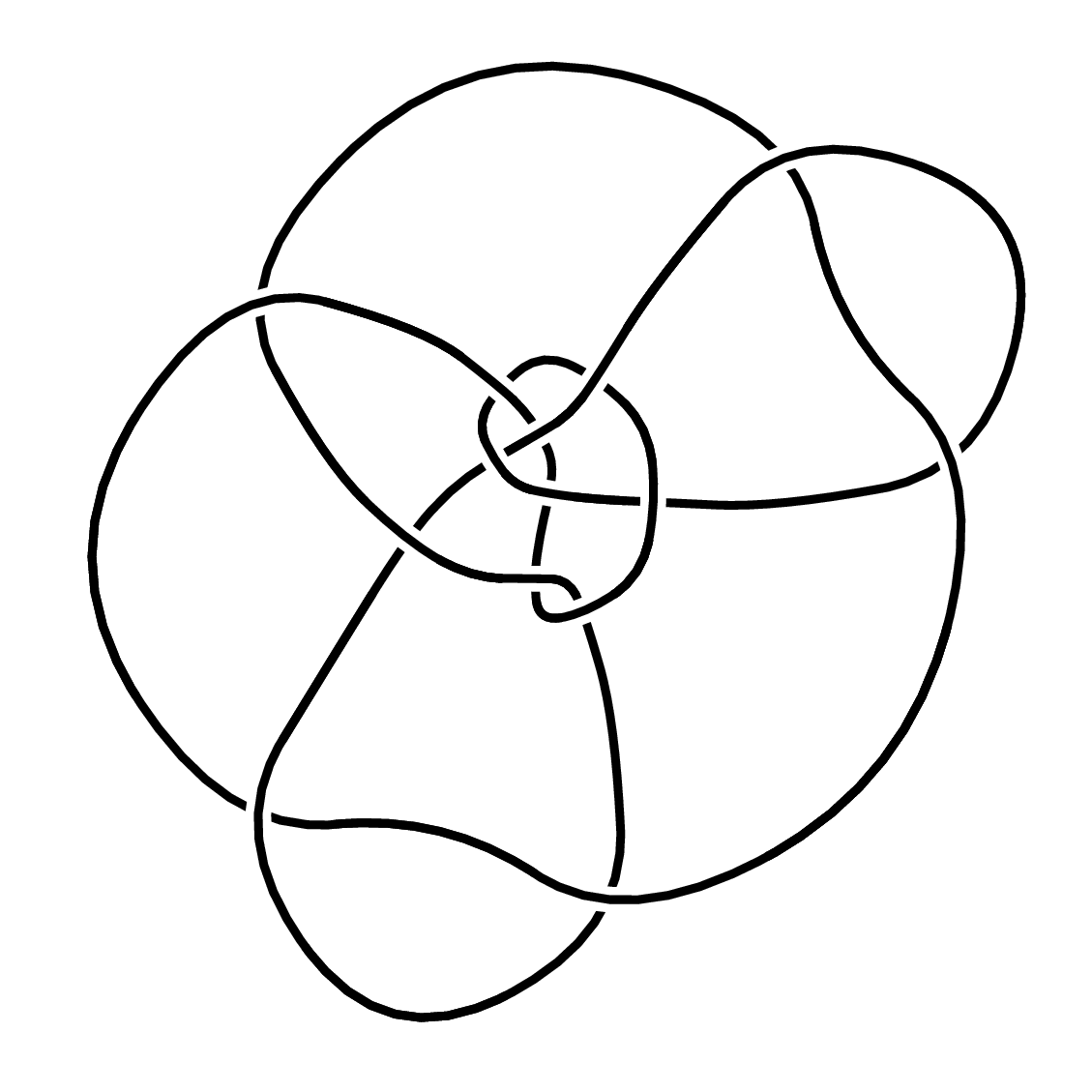}\qquad
\includegraphics[scale=.3]{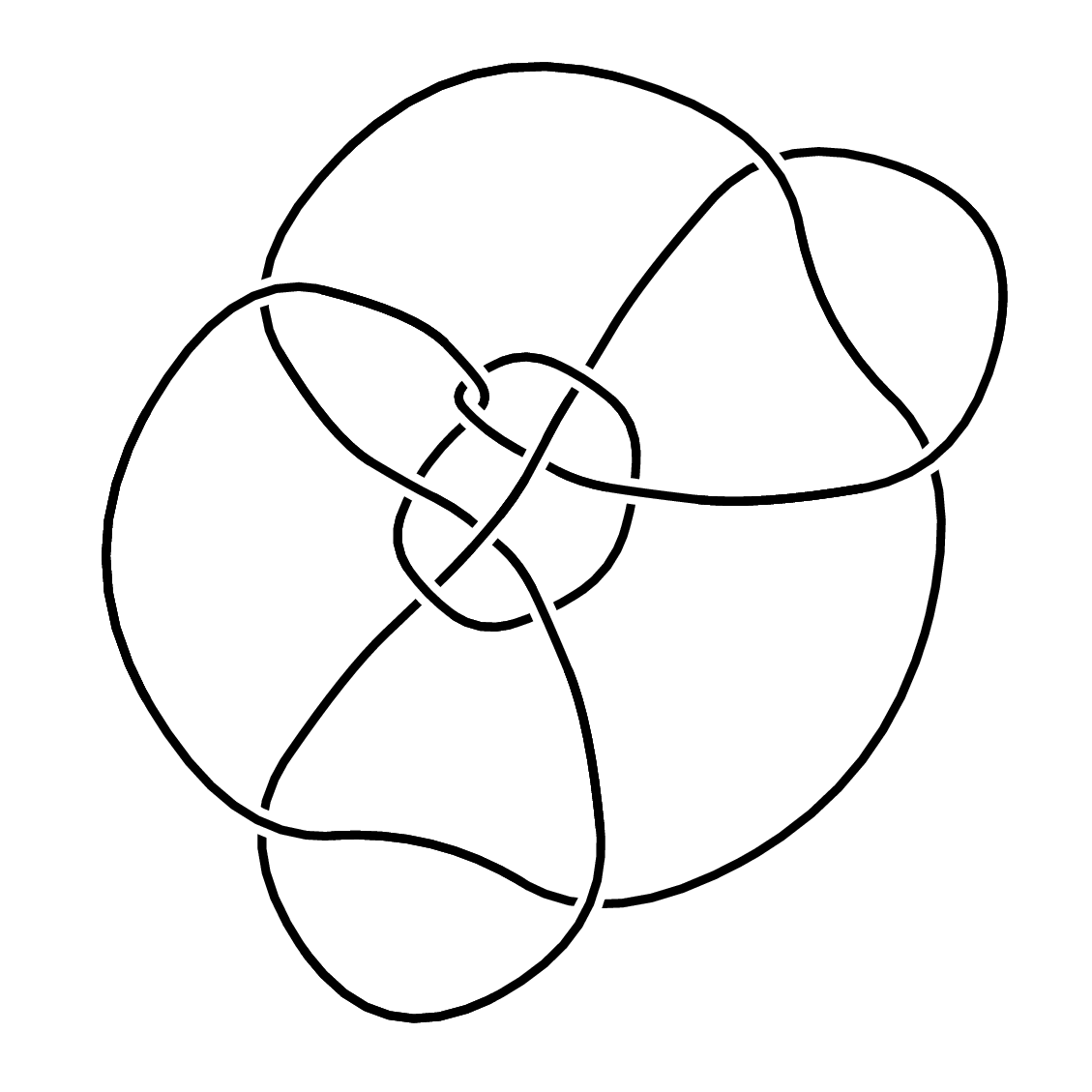}}
\caption{The pair of knots $14^n_{22185}$ (left) and $14^n_{22589}$ (right),
in Knotscape notation.}\label{14npair}
\end{figure}

One interesting open problem about $(0,4)$-mutation is whether this
operation can change the $sl_2$-Khovanov homology introduced
in~\cite{Kh}.   For $(2,0)$-mutation, we settle the analogous question:

\begin{theorem}
  \label{thm.1} The $sl_2$-Khovanov Homology is not invariant under
  $(2,0)$\hyp mutation of knots.  In particular, the pair of
  $(2,0)$-mutant knots in Figure~\ref{14npair} have differing Khovanov
  homologies.
\end{theorem}

For the odd variant of $sl_2$-Khovanov homology, Bloom recently showed
that it is invariant under $(0,4)$-mutation \cite{B}; as a consequence,
the normal $sl_2$\hyp Khovanov homology with mod 2 coefficients is also
invariant.  We do not whether either of these invariants is preserved
by $(2,0)$-mutation.

\begin{question}
  Is the odd $sl_2$-Khovanov homology preserved by genus 2 mutation? 
\end{question}

The $sl_n$-homology introduced by Khovanov and Rozansky~\cite{KR}
cannot be invariant under $(2,0)$-mutation, simply because the Euler
characteristic need not be, since the HOMFLY-PT polynomial can change
under $(2,0)$-mutation.

One final result of this paper is
\begin{proposition}
  \label{prop.1} There exist knots with same colored Jones polynomials
  (for all colors), HOMFLY-PT and Kauffman polynomials, volume and
  signature, but different Khovanov (and reduced Khovanov) homology.
\end{proposition}
The knots from Figure~\ref{14npair} are again examples here, and all
the above claimed properties except for the Khovanov homology are
consequences of the fact that they are $(2,0)$-mutant
(see~Figure~\ref{fig:cabled-mutants}.a).  These same knots were studied by
Stoimenow and Tanaka~\cite{ST1,ST2}, who showed that these knots are not
$(0,4)$-mutants, yet have the same colored Jones polynomials.  (Stoimenow and
Tanaka use notation $14_{41721}$ and $14_{42125}$ for what we denote
$14^n_{22185}$ and $14^n_{22589}$, respectively.)

There are other invariants whose behavior under genus 2 mutation it
would be interesting to understand.  In particular:
 \begin{question}
   Is the Kauffman polynomial invariant under genus 2 mutation?  What
   about the property of having unknotting number one?
\end{question}

Classical Conway $(0,4)$-mutation preserves both these properties
\cite{Lickorish, GordonLuecke2}.  As we discuss in
Section~\ref{sub.kauffman} below, we expect that, in analogy with what
happens with the HOMFLY-PT polynomial, genus 2 mutation should be able
to change the Kauffman polynomial.  \textbf{Addendum:} Morton and
Ryder have confirmed this, showing that the Kauffman polynomial is
\emph{not} invariant under genus 2 mutation \cite{MR2}.

We now detail where the results in Table~\ref{table.main} come from.
The invariance of the hyperbolic volume, or more generally the Gromov
norm, was proven by Ruberman for all types of mutation~\cite{Ru}.  The
statement~\cite[Thm.~1.5]{Ru} requires an additional hypothesis on
$F$, but arguments elsewhere in~\cite{Ru} negate the need for this;
see our discussion of Theorem~\ref{thm.volume} below.  Cooper and
Lickorish proved the invariance of the Alexander polynomial and
generalized signature under a more limited class of $(2,0)$-mutations
than we consider here~\cite{CL}.  This class, which we call handlebody
mutations, turns out to be the main case anyway, and thus it is not
hard to conclude the more general result; see Theorem~\ref{thm.alex}
below.  In the case of the colored Jones polynomials (for a definition
see e.g.~\cite{J,Tu1}), the result essentially follows from
Morton-Traczyk~\cite{MT}, which we modify as Theorem~\ref{thm.cjones}.
In the case of the non-invariance of the HOMFLY-PT polynomials, we
give explicit examples based on the ideas of
Section~\ref{sub.HOMFLY-PT}.

As usual, the presentation of our results does not follow the
historical order by which they were discovered. The project started by
running a computer program of A.Sh. (see~\cite{Sh}) to all knots with
less than or equal to $16$ crossings, taken from Knotscape~\cite{HTh}.
The computer found a single pair of $14$ crossing knots with the same
HOMFLY-PT polynomial, Kauffman polynomial, signature, volume and
different Khovanov Homology, and four pairs of $15$ crossing knots
with same behavior. The knots were isolated, redrawn, and a pattern
was found. Namely, the knots in the above pairs have diagrams that
differ by a so-called cabled mutation (see
Section~\ref{sub.cabledmutation} for a definition). Cabled mutation
can always be achieved by $(2,0)$-mutation.  This, together with a
Kauffman bracket skein theory argument (which we later found in
Morton-Traczyk's work~\cite{MT}) implies that these pairs have
identical colored Jones polynomials, for all colors. At that time, the
numerical equality of the volumes of these pairs was rather
mysterious.  Later on, we found that cabled mutation is a special case
of $(2,0)$-mutation. Ruberman's theorem explained why these pairs have
equal volume.  Once it was observed that Khovanov homology was not
invariant under $(2,0)$-mutation, we asked whether this was true for
other well-known knot invariants, such as the colored Jones
polynomials, the HOMFLY-PT and the Kauffman polynomials.  Once we
realized that the HOMFLY-PT and Kauffman polynomials ought to detect
$(2,0)$-mutation (and even cabled mutation), we tried to find examples
of such knots.

\subsection*{Acknowledgments}
The authors wish to thank I. Agol, D. Bar-Natan and G. Masbaum for
useful conversations; L. Kauffman, J. Przytycki and F. Souza for
organizing an AMS meeting in Snowbird, Utah, and G. Masbaum and P.
Vogel for their hospitality in Paris VII, where the work was
initiated. Finally, we wish to thank the computer team at Georgia Tech
and in particular Lew Lefton and Justin Filoseta for their support in
large scale computations.

\section{The topology of knot mutation}
\label{sec.topologymutation}

This section gives the basic topological lemmas about mutation that we
will need.  In addition to checking that $(2,0)$-mutation of a knot in
$S^3$ makes sense (i.e.~mutating $S^3$ along such a surface always
gives back $S^3$), we will show that one can usually reduce to the
case where the mutation surface has a number of special properties.
Finally, we introduce the notion of cabled mutation for knots in
$S^3$, which is a special type of genus 2 mutation which is easy to
realize diagrammatically.

We begin in the context of general 3-manifolds before specializing to the
case of knots in $S^3$.  From a topological point of view, it is often
best to work with mutation surfaces that are incompressible.  The
following proposition is implicit in~\cite[Sec.~5]{Ru}, and explicit in
a slightly weaker form in~\cite[Lem.~2.2]{Ka2}; one application
below will be to show that mutation makes sense for knots in $S^3$.

\begin{proposition}
  \label{prop.incompressiblemutation} Let $F$ be a closed genus $2$
  surface in a compact orientable 3-manifold $M$. Then either:
\begin{enumerate}
\item
$F$ is incompressible, or
\item
$\Mtau$ can be obtained by mutating along one or two incompressible, 
non-boundary parallel tori, or
\item
$\Mtau \cong M$.
\end{enumerate}
\end{proposition}

\begin{proof}
  The basic idea here is that if $F$ is compressible, then $M^\tau$ is
  homeomorphic to the result of mutating $M$ along any surface
  obtained by compressing $F$.  So suppose $D$ is an embedded
  compressing disc for $F$.  Initially, let us suppose that $\partial D$ is
  a non-separating curve in $F$.  The key property of the
  hyper-elliptic involution $\tau$ is that if $\gamma$ is any
  non-separating simple closed curve in $F$, then $\tau(\gamma)$ is isotopic
  to $\gamma$ with the orientation reversed.  Thus, we can isotope $D$
  so that $\tau(\partial D) = \partial D$, and the restriction of $\tau$ to $\partial D$ is
  a reflection (that is, conjugate to reflecting a circle centered at
  the origin of $\mathbb{R}^2$ about the $x$-axis).
  
  Now perform a surgery of $F$ along $D$ to obtain a surface $T$, which
  consists of the union of $F \setminus N(\partial D)$ with two parallel
  copies of $D$.  Since
  $\partial D$ is non-separating, $T$ is a torus.  There is a natural
  homomorphism $\sigma$ of $T$ which agrees with $\tau$ on $F \setminus N(\partial D)$
  and permutes the two copies of $D$.   We claim that 
  \begin{enumerate}
  \item The involution $\sigma$ is just the elliptic involution of the
    torus shown in Figure~\ref{fig:symmetricsurface}.
    \item $M^\tau \cong M^\sigma$.  
  \end{enumerate}
  The first point is clear, and so turning to the second let us assume
  (for notational simplicity only) that $F$ separates $M$.  Denote by
  $M_1$ and $M_2$ the two pieces of $M$ cut along $F$.  Let $X$ be the
  complement in $M_2$ of a product regular neighborhood $N$ of $D$; we
  can then view our surface $T$ as $\partial X$.  Both $M^\tau$ and $M^\sigma$ can
  be thought of as obtained by gluing together the pieces $M_1$, $X$,
  and $N$.  Moreover, the way that $M_1$ and $X$ are glued is exactly
  the same in both cases, since $\tau$ and $\sigma$ agree on $F \setminus N$; hence
  $M^\tau$ and $M^\sigma$ differ only in how the ball $N$ is attached.
  Since there is a unique way of attaching a 3-ball to a 2-sphere up
  to homeomorphism, we have $M^\tau \cong M^\sigma$ as claimed.   (You can also
  see the homeomorphism of $N$ needed to build the map $M^\tau \to
  M^\sigma$ directly --- thinking of $N$ as a pancake, just flip it
  over.)   
  
  Thus in the case that $\partial D$ is non-separating, we have shown that
  $M^\tau$ is homeomorphic to a mutant of $M$ along a torus $T$.  If
  $\partial D$ is separating, then the picture is essentially the same.  In
  this case, we can isotope $\partial D$ so that $\tau$ fixes it pointwise.
  Proceed as above, the only difference being that now surgering
  $F$ along $D$ results in a disconnected surface consisting of two
  tori.  Thus in either case, $M^\tau$ is homeomorphic to the result
  of mutating $M$ along either one or two tori.
  
  So to complete the proof of the proposition, we just need to show
  that if $T$ is a torus in $M$ with elliptic involution $\sigma$,
  then either
  \begin{enumerate}
    \item $T$ is incompressible and not boundary parallel. 
    \item $M^\sigma \cong M$.  
  \end{enumerate}
  If $T$ were boundary parallel, then mutating along it doesn't change
  the topology since the gluing map $\sigma$ extends over the product
  region bounded by $T$ and a component of $\partial M$.  If $T$ is
  compressible, then arguing as above we see that $M^\sigma$ is
  homeomorphic to the result of mutating along a 2-sphere $S$ in $M$,
  where the gluing map $\phi$ is just rotation of $S$ about some axis
  through angle $\pi$; since $\phi$ is isotopic to the identity, we
  have that $M \cong M^\phi \cong M^\sigma$, as desired.  
\end{proof}

\begin{remark}\label{remark.boundaryslopes}
  Later, we will apply this proposition to a manifold $M$ where $\partial M$
  is a torus, and need the following fact.  As setup, note
  that since $F$ is closed, there is a canonical identification of
  $\partial M$ with $\partial M^\tau$.  The observation is that if we end up in case
  (3) where $M^\tau \cong M$, then the proof shows that there is a
  homeomorphism $f :  M \to M^\tau$ where the restriction of $f$ to
  $\partial M$ is either the identity or the elliptic involution.  (The later
  happens when part of $F$ compresses to something parallel to the
  boundary torus.)
\end{remark}

\begin{remark}\label{rem.othersurfaces}
  While Proposition~\ref{prop.incompressiblemutation} nominally
  concerns only genus 2 mutation, there are analogous statements for
  any of the symmetric surfaces, which follow from the same proof.
\end{remark}

Ruberman proved that if $M$ is hyperbolic, and $F$ any symmetric
surface in $M$, then $M^\tau$ is also hyperbolic and, moreover, $M$ and
$M^\tau$ have the same volume.  This is stated in~\cite[Thm.~1.3]{Ru}
with the additional hypothesis that $F$ is incompressible.  However, as
he observed in Section~5 of that same paper, this hypothesis can be
dropped by appealing to Proposition~\ref{prop.incompressiblemutation}
and Remark~\ref{rem.othersurfaces}.   Similarly, one has:

\begin{theorem}[\cite{Ru}]
  \label{thm.volume} Let $M$ be a orientable 3-manifold, whose boundary,
  if any, consists of tori.  Then the result of mutating $M$ along any
  symmetric surface has the same Gromov norm as $M$ itself.  
\end{theorem}

In the context of knots in $S^3$ that we consider below, we will be
dealing with manifolds where $\partial M$ is a single torus.  In this case,
Ruberman~\cite[Sec.~5]{Ru} and Tillmann~\cite[Rem.~1.3]{Ti1} observed
that all of the types of mutations pictured in
Figure~\ref{fig:symmetricsurface} can be reduced to a sequence of
genus 2 mutations, provided the mutation surface is separating.  
\begin{lemma}[\cite{Ru, Ti1}]
  \label{lem.20suffice} Suppose $M$ is a compact orientable 3-manifold
  whose boundary is a single torus.  Let $F$ be one of the symmetric
  surfaces depicted in Figure~\ref{fig:symmetricsurface}.  Provided $F$
  is separating, mutation along $F$ can always be accomplished by a
  composition of at most two $(2,0)$-mutations.
\end{lemma}
The idea they used to prove this lemma is to tube copies of $F$ along
$\partial F$ to build a closed genus 2 surface $S$.  Mutating along $S$ is
the same as doing a certain mutation along the original surface
$F$, for reasons similar to the proof of
Proposition~\ref{prop.incompressiblemutation}.  In the case where $F$
is a 4-punctured sphere, it may not be possible that the desired
involution $\tau_i$ can be directly induced by mutation along a tubed
surface $S$; however, in this case the needed  mutation can be realized
by mutating along the possible choices for $S$ in succession.  

\subsection{Genus 2 mutation of knots in $S^3$}\label{subsec.knotmut}

Suppose that $F$ is a closed genus 2 surface in $S^3$.  As $S^3$ is
simply connected, the Loop Theorem implies that $F$, as well as any
torus in $S^3$, is compressible.  Therefore, the trichotomy of
Proposition~\ref{prop.incompressiblemutation} forces $(S^3)^\tau$, the
result of mutation along $F$, to again be homeomorphic to $S^3$.  Thus
if $K$ is a knot in $S^3$ disjoint from $F$, then we can consider the
resulting knot $K^\tau$ in $(S^3)^\tau \cong S^3$, which we call the mutant of $K$ along $F$. 

When the surface $F$ bounds a genus 2 handlebody $H$ in $S^3$, then
the mutation operation is particularly simple to describe, since the
hyper-elliptic involution $\tau$ extends to give a self-homeomorphism
of $H$.  When the knot $K$ is contained in $H$, we say that $K^\tau$ is
obtained from $K$ by $(2,0)$-{\em handlebody mutation}.  (If instead
$K$ is in the complement of $H$, then $K^\tau \cong K$.)  Such
$(2,0)$-handlebody mutation was studied by Cooper-Lickorish~\cite{CL},
who were interested in how it affected the Alexander polynomial.

The next proposition shows that $(2,0)$\hyp handlebody mutation is
actually the main interesting case of genus 2 mutation, the only other
case being $(1,0)$\hyp {\em handlebody mutation}, which is defined
analogously.

\medskip
\noindent
\textbf{This proposition is false, and hence following corollary has
  not been established.  See the attached erratum for correct proofs
  of the results depending on Proposition 2.7 and Collorary 2.8.}

\begin{proposition}
  \label{prop.handlebodymutation} Let $K$ be a knot in $S^3$ which is
  disjoint from a genus $2$ surface $F$. Then either:
\begin{itemize}
\item
$\Ktau$ is obtained from $K$ by $(2,0)$-handlebody mutation, or
\item
$\Ktau$ is obtained from $K$ by one or two $(1,0)$-handlebody mutations, or
\item
$\Ktau \cong K$.
\end{itemize}
\end{proposition}

\begin{proof}
  Let $M = S^3 \setminus N(K)$ be the exterior of $K$.  Applying
  Proposition~\ref{prop.incompressiblemutation} to $F$ thought of as
  a surface in $M$, we have three cases.
 
  First, $F$ may be incompressible in $M$; in this case, we claim this
  is actually a $(2,0)$-handlebody mutation.  Let $X$ and $Y$ be the two
  pieces of $S^3$ cut along $F$, and suppose that $K$ lies in $X$.
  Since $F$ is incompressible in $M$, it is also incompressible as the
  boundary of $Y$.  Thus any compressing disc for $F$ in $S^3$ lies in
  $X$.  Pick two such compressing discs, whose boundaries are disjoint
  non-parallel non-separating curves in $F$ (by Dehn's Lemma, every
  embedded curve in $F$ bounds a compressing disc as $\pi_1(S^3) = 1$).
  If we compress $F$ along both these discs, we get a sphere which
  bounds a ball on \emph{both} sides.  This shows $X$ is handlebody.
  
  Second, suppose mutation along $F$ in $M$ can be achieved by one or
  two mutations along incompressible tori.  The argument just given
  shows that those are $(1,0)$-handlebody mutations.
  
  Finally, suppose that we are in the final case where $M^\tau \cong M$.  This
  shows that the complements of $K^\tau$ and $K$ are the same, but we
  need to show that the knots themselves are the same.  Of course,
  knots are determined by their complements~\cite{GordonLuecke},
  but we now give an elementary argument.  We can reconstruct $K$ from
  $M$ if we just mark the loop on $\partial M$ which is the meridian for $K$, and the
  same for $K^\tau$ and $M^\tau$.  By Remark~\ref{remark.boundaryslopes},
  the homeomorphism of $M^\tau \to M$ takes the meridian to the meridian,
  establishing $K^\tau \cong K$ as desired.
\end{proof}

A $(1,0)$-handlebody mutation may be realized by a $(2,0)$-handlebody mutation
simply by adding a nugatory handle. Thus,

\medskip
\noindent
\textbf{This corollary has not been established; see the comment
  above Proposition 2.7.}

\begin{corollary}
\label{cor.handlebodymutation}
Any knot invariant which
does not change under $(2,0)$\hyp handlebody mutation, does not change under
$(2,0)$-mutation.
\end{corollary} 

Using this, we can generalize~\cite{CL} to:

\begin{theorem}\label{thm.alex}
  The Alexander polynomial and the generalized signatures of a knot in
  $S^3$ do not change under $(2,0)$-mutation.
\end{theorem}

\begin{proof}  
  In~\cite[Cor.8]{CL} Cooper-Lickorish prove that these invariants do
  not change under $(2,0)$-handlebody mutation. The result thus follows
  from Corollary \ref{cor.handlebodymutation}.
\end{proof}

\subsection{Cabled mutation}
\label{sub.cabledmutation}

In this short section, we introduce the notion of cabled mutation,
which is a special form of genus 2 mutation which we will use to
construct examples where the HOMFLY-PT polynomial
changes under mutation.  

Consider a framed 2-2 tangle $T$ in a ball, that is, a ball containing
two disjoint properly embedded arcs (the \emph{strings}), where each
arc has a preferred framing.  If $T$ were part of a knot, then we
could do $(0,4)$-mutation on it using one of the three involutions
pictured in Figure~\ref{fig:symmetricsurface}.  Let $\tau$ be one of
these involutions which is string-preserving, that is, exchanges one
of the endpoints of a fixed arc with the other.  Let $\Ttau$ denote
the image of $T$ under the involution.  Given natural numbers $n,m \geq
1$, let $T(n,m)$ (resp.~$\Ttau(n,m)$) denote the tangle obtained by
taking a $n$ and $m$ parallel of the strings of $T$ (resp.~$\Ttau$).

\begin{definition}
  \label{def.cabledmutation} Connected cabled mutation (or simply,
  cabled mutation) is the result of replacing $T(n,m)$ by $\Ttau(n,m)$ in
  some planar diagram of a knot in $S^3$.
\end{definition}

When $n=m=1$, cabled mutation is just usual $(0,4)$-mutation.  One
motivation for studying this notion is that $(0,4)$-mutation followed
by connected cabling can be often be achieved by a connected cabled
mutation.

Our next lemma discusses the relation between cabled mutation and
genus 2 mutation.

\begin{lemma}
\label{lem.cabledmutationgenus2}
Cabled mutation is a
special form of genus 2 mutation. 
\end{lemma}

\begin{proof} 
Starting with the boundary of the
tangle $T$ we can attach two tubes inside it, containing the strands
of $T(n,m)$, to produce a closed genus 2 surface $F$.  The cabled
mutation on $T(n,m)$ can then be achieved by cutting along $F$ and
regluing; because the original involution on $T$ is string preserving,
the map we reglue $F$ by is the hyper-elliptic involution $\tau$
pictured in Figure~\ref{fig:symmetricsurface}.  (If $\tau$ was not
strand preserving, then the regluing map for $F$ is some other
involution and this is not a mutation.)
\end{proof}

\section{Behavior of quantum invariants under mutation}
\label{sec.TQFT}

As mentioned in the introduction, many knot invariants are preserved
under Conway $(0,4)$-mutation.  Such invariants include the HOMFLY-PT
(and, hence, Jones and Alexander) and Kauffman polynomials, see for
example~\cite{Lickorish, LL,MC,MT, CL}. In this section we analyze the
behavior of several quantum invariants under $(2,0)$-mutation.

\subsection{Invariance of the Jones polynomials under $(2,0)$-mutation}
\label{sub.jones}

Morton and Traczyk showed that the colored Jones polynomials are
invariant under Conway mutation \cite{MT}.  As we now describe, their
approach easily generalizes:
\begin{theorem}
\label{thm.cjones}
The colored Jones polynomials of a knot are invariant under $(2,0)$-mutation
for all colors.
\end{theorem}

\begin{proof}
  The theorem follows from the fact that the colored Jones polynomial
  can be defined via the Kauffman bracket skein theory, in the style
  of topological quantum field theory, see~\cite{Kf}.  By
  Corollary~\ref{cor.handlebodymutation} it suffices to consider genus
  2 \emph{handlebody} mutation.  
  \begin{figure}[htb]
    \centerline{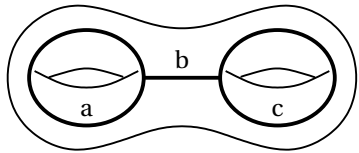}
    \caption{Basis of the Kauffman skein module of a closed genus $2$ surface.}
    \label{fig:skein-gens}
  \end{figure}

  The Kauffman bracket skein module of a genus $2$ handlebody has a
  basis that consists of all the colored trivalent graphs $G(a,b,c)$,
  where $a$, $b$, and $c$ are nonnegative integers with $c \leq 2
  \min\{a,b\}$ (see Figure~\ref{fig:skein-gens}).  Indeed, a genus 2
  handlebody is diffeomorphic to a $(\text{twice punctured disk}) \times
  I$, and a basis for the Kauffman bracket of the latter is given in
  \cite[Cor.~4.4]{PS}.  Since this basis is clearly invariant under
  $\tau$, it implies that the colored Jones polynomials are invariant
  under $(2,0)$-handlebody mutation, proving the theorem.
\end{proof}

Combining Theorem~\ref{thm.cjones} with the Melvin-Morton-Rozansky
Conjecture (settled in~\cite{B-NG}) gives an alternate proof of
Theorem~\ref{thm.alex}, namely that the Alexander polynomial of a knot
is invariant under $(2,0)$-mutation.

\subsection{Non-invariance of HOMFLY-PT under $(2,0)$-mutation}
\label{sub.HOMFLY-PT}

It is not hard to see that the HOMFLY-PT and Kauffman polynomials are
invariant under $(0,4)$\hyp mutation \cite{Lickorish}. This follows from
the fact that the corresponding skein modules of a 3-ball with 4 marked
points on the boundary have a basis consisting of the following three
diagrams that are invariant under the involution in question:

$$ 
\includegraphics[width=0.5in]{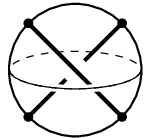}\qquad\qquad
\includegraphics[width=0.5in]{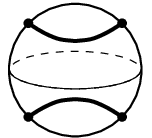}\qquad\qquad
\includegraphics[width=0.5in]{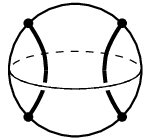}
$$

\begin{figure}[ht]
\centerline{\includegraphics[width=4.0in,trim=0 250 0 0]{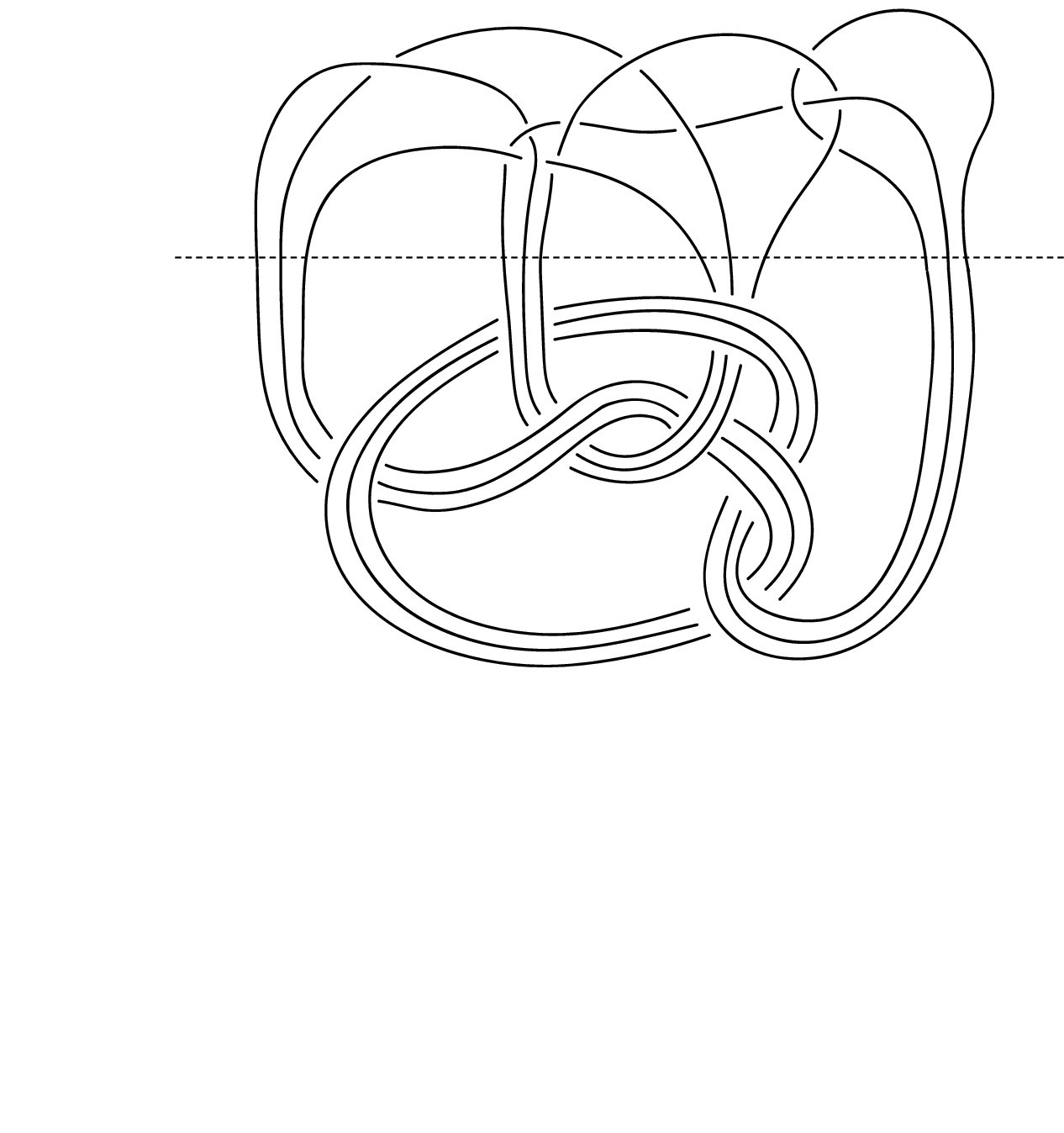}}
\caption{The knot $K_{75}$. 
It and its cabled mutant $K_{75}^\tau$ have different
HOMFLY-PT polynomials}
\label{fig:cabled-HOMFLY-PT}
\end{figure}

In contrast, genus 2 mutation can change the HOMFLY-PT polynomial.  In
particular, we found a 75 crossing knot $K_{75}$ which has a cabled
mutant with differing HOMFLY-PT polynomials.  This knot is depicted in
Figure~\ref{fig:cabled-HOMFLY-PT}.  As you can see, $K_{75}$ contains
a $(3,3)$-cabled tangle which is the region below the horizontal line;
let $K_{75}^\tau$ be the cabled mutant of $K_{75}$ with respect to a
string-preserving involution $\tau$ of this tangle.    

Direct computation with the Ewing-Millett computer
program implemented in Knotscape shows that $K_{75}$ and $K_{75}^\tau$ have
different HOMFLY-PT polynomials. Coefficients of these polynomials are given
in Tables~\ref{tbl:HOMFLY-PT1} and \ref{tbl:HOMFLY-PT2} on pages~\pageref{tbl:HOMFLY-PT1} and \pageref{tbl:HOMFLY-PT2} (with zero
entries omitted). For example, the coefficient of the monomial $m^2 l^{-2}$ is
$56$ in both polynomials. On the other hand, the coefficients of $m^4 l^{-2}$
are $-953$ for $K_{75}$ and $-964$ for $K_{75}^\tau$.

%
%
%

Here is a heuristic reason why the HOMFLY-PT polynomial is not
invariant under $(2,0)$-mutation, which explains how we came across
our pair of $75$ crossing knots.  First, it was already known that
there are $(2,0)$-mutant \emph{links} with different HOMFLY-PT
polynomials \cite{CL}.  In particular, start with the
Kinoshita-Terasaka and Conway knots which are a famous pair of 11
crossing knots which differ by $(0,4)$-mutation.  Morton and Traczyk
showed (see \cite{MC}) that taking a certain disconnected $3$-cable of
each of these knots gives a pair of links with differing HOMFLY-PT
polynomials; this gives a pair of cabled-mutant links with distinct
HOMFLY-PT polynomials.  (In contrast, Lickorish-Lipson
showed~\cite{LL} that the HOMFLY-PT polynomial of 2-cables of mutant
knots are always equal.)  This suggests that we should have a good
chance of getting a pair of connected cabled mutant knots with
distinct HOMFLY-PT polynomials by the following procedure: take as a
pattern tangle the one that appears in the Kinoshita\hyp Terasaka and
Conway pair, cable each of its components $3$ times, and close it up
to a knot in some fairly arbitrary way.  This is exactly how we found
the pair of knots with $75$ crossings.

\subsection{Expected non-invariance of the Kauffman polynomial}
\label{sub.kauffman}

The heuristic reasons for the non-invariance of the HOMFLY-PT
polynomial under $(2,0)$\hyp mutation applies equally well in the case of
the Kauffman polynomial.  For this reason, we expect that the Kauffman
polynomial is not invariant under $(2,0)$-mutation.  To show this, it
suffices to present a pair of cabled mutant knots with different
Kauffman polynomials.  However, the available computer programs for
computing the Kauffman polynomial do not work well with knots with
more than 50 crossings, and this has prevented us from examining any
interesting examples.  \textbf{Addendum:} Morton and
Ryder have now succeeded in showing that the Kauffman polynomial is
\emph{not} invariant under genus 2 mutation \cite{MR2}.

\subsection{Proof of Proposition~\ref{prop.1}}
\label{sub.prop1}

\begin{figure}[htb]
\scalebox{0.95}{\centerline{\qquad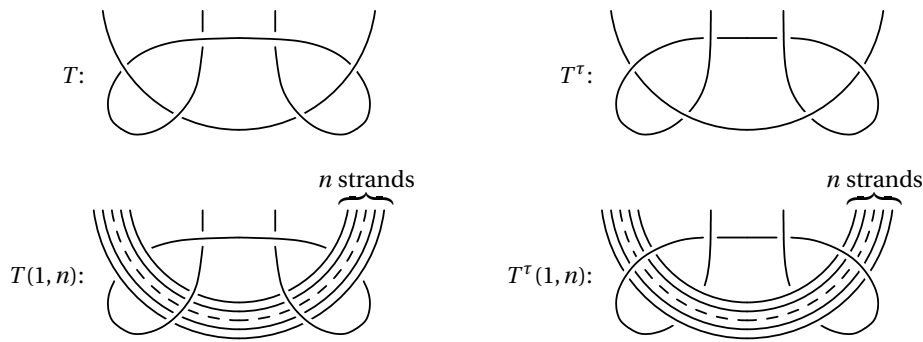}}
\caption{Cabling of a tangle and its mutant.}\label{fig:tangles-cabled}
\end{figure}

Now we show there exist knots with the same colored Jones, HOMFLY-PT,
and Kauffman polynomials, the same volume and signature, but different
Khovanov homology.  Consider the tangles $T$ and $\Ttau$ from
Figure~\ref{fig:tangles-cabled}. Denote by $T(1,n)$ and $\Ttau(1,n)$
their $(1,n)$-cables, respectively (for some fixed $n$). Let $K$ and
$\Ktau$ be two knots that differ by replacement of $T(1,n)$ with
$\Ttau(1,n)$.  In particular, $K$ and $\Ktau$ are connected cabled
mutants and, thus, $(2,0)$-mutant.  Theorems \ref{thm.volume} and
\ref{thm.cjones} thus imply that $K$ and $\Ktau$ have equal colored
Jones polynomials and volume.  A priori, $K$ and $\Ktau$ could have
different HOMFLY-PT and Kauffman polynomials.  However, an elementary
computation in the respective skein theories imply that $K$ and
$\Ktau$ also have equal HOMFLY-PT and Kauffman polynomials.

\begin{figure}[htb]
\scalebox{0.75}{\centerline{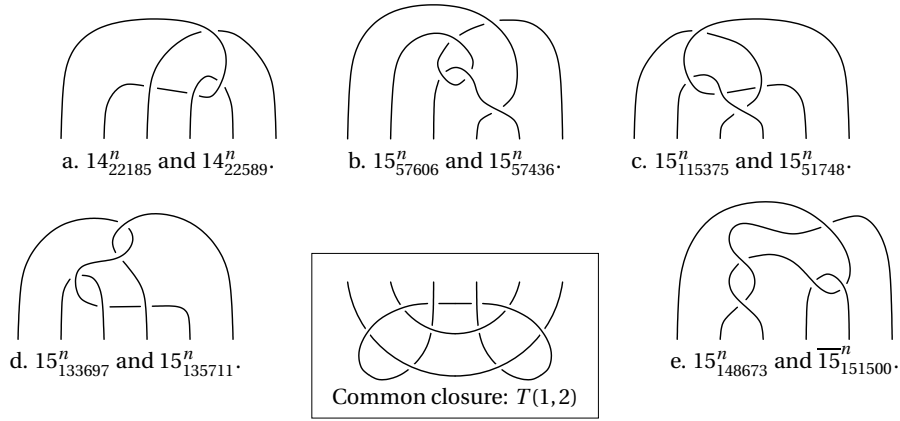}}
\caption{Five pairs of cabled mutant knots with at most 15 crossings that have
different Khovanov homology. They are closures of the $T(1,2)$ tangle.}
\label{fig:cabled-mutants}
\end{figure}

When $n=2$, let us choose the closure of $T(1,2)$ in one of the ways from
Figure~\ref{fig:cabled-mutants} to obtain five pairs of knots. In Knotscape
notation~\cite{HTh}, these five pairs are $(14^n_{22185},14^n_{22589})$,
$(15^n_{57606},15^n_{57436})$,
$(15^n_{115375},15^n_{51748})$,
$(15^n_{133697},15^n_{135711})$, and
$(15^n_{148673},\overline{15}^n_{151500})$,
where the bar above the number of crossings means the mirror
image of the corresponding knot. Computer calculations with {\tt
KhoHo}~\cite{Sh} show that knots from these pairs have different Khovanov
Homology (see Section~\ref{sub.data}).

\subsection{Knots with few crossings}
\label{sub.data}

We say that two knots are {\em almost mutant} if they have the same
HOMFLY-PT and Kauffman polynomials, signature, and hyperbolic volume.
This is an equivalence relation. Note that mutant knots are almost
mutant.

We can partition the set of knots with a bounded number of crossings
according to the equivalence relation of being almost mutant.  We
worked out these equivalence classes for all knots with at most $16$
crossings. As it turns out, almost mutant knots with at most 16
crossings always have the same number or crossings. As a consequence,
two such knots are either both alternating or both non-alternating.
This follows from the fact that the span of the Jones polynomial of a
knot equals the number of crossings for this knot if and only if the
knot is alternating.  For non-alternating knots,
Table~\ref{table:AM_classes} lists the number of such equivalence
classes of a given size.  We restrict the table to non-alternating
knots only because we are interested primarily in the possibilities for
the Khovanov homology of almost mutant pairs; for alternating knots,
the Khovanov homology (at least the free part thereof) is completely
determined by their Jones polynomials and signature~\cite{L}.

The number of knots in Table~\ref{table:AM_classes} is taken from Knotscape,
which does not distinguish between mirror images. Therefore, we considered
each knot twice: the knot itself and its mirror image. The number of
amphicheiral knots can be found in~\cite{HThW}. The notation
$a_1:n_1, a_2:n_2,\,\dots\,,a_k:n_k$ means that there are $n_j$ equivalence
classes of size $a_j$ for $j=1,2,\,\dots\,,k$.

\begin{table}[t]
\begin{center}{\small
\begin{tabular}{|p{0.85cm}|p{1.24cm}|p{1.35cm}|p{0.95cm}|p{6.2cm}|} \hline
\centering num. cross. &\centering num. knots & \centering counting mirrors & amph. knots &\centering size and number \\  of almost mutant classes \tabularnewline
\hline\hline
\hfill $\leq13$ &    \hfill 6236 &  \hfill 12468 &  \hfill 4 & 2:\,1028, 3:\,54, 4:\,42, 6:\,2 \\ \hline
    \hfill $14$ &   \hfill 27436 &   \hfill 54821 & \hfill 51 & 2:\,5349, 3:\,298, 4:\,359, 6:\,30, 8:\,10 \\ \hline
    \hfill $15$ &  \hfill 168030 &  \hfill 336059 &  \hfill 1 & 2:\,35423, 3:\,1368, 4:\,4088, 6:\,290, 8:\,136 \\ \hline                         
    \hfill $16$ & \hfill 1008906 & \hfill 2017322 &\hfill 490 &  2:\,212351, 3:\,6612, 4:\,33156, 6:\,2159, 7:\,20, 8:\,2229, 9:\,4, 10:\,8, 12:\,201, 16:\,22, 20:\,2 \\
\hline
\end{tabular}
}\end{center}
\caption{Sizes and numbers of almost mutant classes of non-alternating knots}
\label{table:AM_classes}
\end{table}

It is remarkable that very few almost mutant knots have different Khovanov
homology. There are only 5 pairs (10 if counted with mirror images) of such
knots with at most 15 crossings. They are exactly the 5 cabled mutant pairs
from Section~\ref{sub.prop1} (see Figure~\ref{fig:cabled-mutants})! We list
values of various knots invariants for these knots below.

There are 27 pairs (54 with mirrors) of almost mutant knots with 16 crossings
that have different Khovanov homology. Many of these pairs consist of cabled
mutant knots, but we could not verify them all. The pairs are:
\begin{center}
\begin{tabular}{ccc}
($16^n_{257474}$, $16^n_{293658}$) & ($16^n_{258027}$, $16^n_{380926}$) & ($16^n_{258035}$, $16^n_{359938}$) \\ 
($16^n_{261803}$, $16^n_{300395}$) & ($16^n_{262535}$, $16^n_{300387}$) & ($16^n_{306846}$, $16^n_{307597}$) \\ 
($16^n_{332130}$, $16^n_{707045}$) & ($16^n_{337388}$, $16^n_{697474}$) & ($16^n_{472161}$, $16^n_{635329}$) \\
($16^n_{564024}$, $16^n_{564036}$) & ($16^n_{564059}$, $16^n_{564068}$) & ($16^n_{789164}$, $16^n_{797712}$) \\
($16^n_{789206}$, $16^n_{797688}$) & ($16^n_{809314}$, $\overline{16}^n_{850490}$) & ($16^n_{809334}$, $\overline{16}^n_{850512}$) \\ 
($16^n_{812818}$, $16^n_{850972}$) & ($16^n_{820956}$, $16^n_{820968}$) &($16^n_{822219}$, $\overline{16}^n_{822229}$) \\
($16^n_{878609}$, $\overline{16}^n_{944604}$) & ($16^n_{884231}$, $16^n_{884268}$) & ($16^n_{885298}$, $16^n_{885312}$) \\
($16^n_{885305}$, $16^n_{885319}$) & ($16^n_{885467}$, $16^n_{885968}$) & ($16^n_{890470}$, $\overline{16}^n_{944600}$) \\
($16^n_{937845}$, $\overline{16}^n_{947575}$) & ($16^n_{939163}$, $16^n_{945493}$) & ($16^n_{943082}$, $16^n_{943119}$) \\
\end{tabular}
\end{center}

We used Knotscape~\cite{HTh} to list all almost mutant knots with at most 16
crossings. Khovanov homology was computed using {\tt KhoHo}~\cite{Sh} for all
knots with at most 15 crossings and  {\tt JavaKh}~\cite{B-NGr} for
non-alternating knots with 16 crossings. It is worth noticing that Knotscape
only computes hyperbolic volume with the precision of 12 significant digits.
We used program Snap~\cite{Snap} to compute the volume with the precision of
180 significant digits to verify our data. As it turns out, there are no knots
with at most 16 crossings that have non-zero difference in hyperbolic volumes
that is less than $10^{-13}$.  Only 132 pairs of knots have difference in
volumes less than $10^{-9}$ and, hence, are considered as having the same 
volume by Knotscape. None of these pairs are almost mutants.

To end this section we list values of some quantum and hyperbolic
invariants for the almost mutant knots with at most 15 crossings that have
different Khovanov homology. They were
computed using Knotscape~\cite{HTh} and {\tt KhoHo}~\cite{Sh}.
HOMFLY-PT and Kauffman polynomials are given by the tables of their
coefficients. Our notation for Khovanov homology is borrowed
from~\cite{B-N2}. An expression $a^i_j$ in the ``ranks'' string means
that the multiplicity of $\BZ$ in the Khovanov homology group with
homological grading $i$ and $q$-grading $j$ is $a$. Negative grading
is shown with underlined numbers. A similar convention is used for
$2$-torsion as well (this is the only torsion that appears for these
knots). In this case, $a$ is the multiplicity of $\BZ_2$. For example,
the homology group of $14^n_{22185}$ with homological grading $0$ and
$q$-grading $(-1)$ is $\BZ^2\oplus\BZ_2^2$.

\vfill 

\expandafter\def\csname Alexander_K14n_22185\endcsname{%
\hbox{$1$}\allowbreak\thinspace
}

\expandafter\def\csname Jones_K14n_22185\endcsname{%
\hbox{$-\,t^{-6}$}\allowbreak\thinspace\hbox{$+\,t^{-5}$}\allowbreak\thinspace\hbox{$+\,t^{-2}$}\allowbreak\thinspace\hbox{$-\,t^{-1}$}\allowbreak\thinspace\hbox{$+\,2$}\allowbreak\thinspace\hbox{$-\,t$}\allowbreak\thinspace\hbox{$-\,t^{4}$}\allowbreak\thinspace\hbox{$+\,t^{5}$}\allowbreak\thinspace
}

\expandafter\def\csname HOMFLYPT_K14n_22185\endcsname{%
$\bgroup\arraycolsep=3pt\begin{array}{|r|rrrrr|}\hline
\vrule width 0pt depth 0pt height 10pt
&\mathbf{l^{-4}}&\mathbf{l^{-2}}&\mathbf{1}&\mathbf{l^{2}}&\mathbf{l^{4}}\\ \hline
\mathbf{       1}&  -3&   8&  -5&    &   1\\
\mathbf{ m^{  2}}&  -4&  14& -11&    &   1\\
\mathbf{ m^{  4}}&  -1&   7&  -6&    &    \\
\mathbf{ m^{  6}}&    &   1&  -1&    &    \\
\hline
\end{array}\egroup$
}

\expandafter\def\csname Kauffman_K14n_22185\endcsname{%
$\bgroup\arraycolsep=1pt\begin{array}{|@{\hspace{3pt}}r@{\hspace{3pt}}|@{\hspace{3pt}}rrrrrrrrrr@{\hspace{3pt}}|}\hline
\vrule width 0pt depth 0pt height 10pt
&\mathbf{a^{-4}}&\mathbf{a^{-3}}&\mathbf{a^{-2}}&\mathbf{a^{-1}}&\mathbf{1}&\mathbf{a}&\mathbf{a^{2}}&\mathbf{a^{3}}&\mathbf{a^{4}}&\mathbf{a^{5}}\\ \hline
\mathbf{       1}&   1&    &    &    &  -5&    &  -8&    &  -3&    \\
\mathbf{ z^{  2}}&    &   5&    &  11&    &  11&    &   9&    &   4\\
\mathbf{ z^{  4}}&  -8&    & -10&    &  15&    &  32&    &  15&    \\
\mathbf{ z^{  6}}&    & -18&    & -30&    & -22&    & -27&    & -17\\
\mathbf{ z^{  8}}&  14&    &  18&    & -15&    & -51&    & -32&    \\
\mathbf{ z^{ 10}}&    &  20&    &  27&    &  13&    &  26&    &  20\\
\mathbf{ z^{ 12}}&  -7&    &  -8&    &   7&    &  35&    &  27&    \\
\mathbf{ z^{ 14}}&    &  -8&    &  -9&    &  -2&    &  -9&    &  -8\\
\mathbf{ z^{ 16}}&   1&    &   1&    &  -1&    & -10&    &  -9&    \\
\mathbf{ z^{ 18}}&    &   1&    &   1&    &    &    &   1&    &   1\\
\mathbf{ z^{ 20}}&    &    &    &    &    &    &   1&    &   1&    \\
\hline
\end{array}\egroup$
}

\expandafter\def\csname HVolume_K14n_22185\endcsname{%
$8.878159662$
}

\expandafter\def\csname Alexander_K14n_22589\endcsname{%
\hbox{$1$}\allowbreak\thinspace
}

\expandafter\def\csname Jones_K14n_22589\endcsname{%
\hbox{$-\,t^{-6}$}\allowbreak\thinspace\hbox{$+\,t^{-5}$}\allowbreak\thinspace\hbox{$+\,t^{-2}$}\allowbreak\thinspace\hbox{$-\,t^{-1}$}\allowbreak\thinspace\hbox{$+\,2$}\allowbreak\thinspace\hbox{$-\,t$}\allowbreak\thinspace\hbox{$-\,t^{4}$}\allowbreak\thinspace\hbox{$+\,t^{5}$}\allowbreak\thinspace
}

\expandafter\def\csname HOMFLYPT_K14n_22589\endcsname{%
$\bgroup\arraycolsep=3pt\begin{array}{|r|rrrrr|}\hline
\vrule width 0pt depth 0pt height 10pt
&\mathbf{l^{-4}}&\mathbf{l^{-2}}&\mathbf{1}&\mathbf{l^{2}}&\mathbf{l^{4}}\\ \hline
\mathbf{       1}&  -3&   8&  -5&    &   1\\
\mathbf{ m^{  2}}&  -4&  14& -11&    &   1\\
\mathbf{ m^{  4}}&  -1&   7&  -6&    &    \\
\mathbf{ m^{  6}}&    &   1&  -1&    &    \\
\hline
\end{array}\egroup$
}

\expandafter\def\csname Kauffman_K14n_22589\endcsname{%
$\bgroup\arraycolsep=1pt\begin{array}{|@{\hspace{3pt}}r@{\hspace{3pt}}|@{\hspace{3pt}}rrrrrrrrrr@{\hspace{3pt}}|}\hline
\vrule width 0pt depth 0pt height 10pt
&\mathbf{a^{-4}}&\mathbf{a^{-3}}&\mathbf{a^{-2}}&\mathbf{a^{-1}}&\mathbf{1}&\mathbf{a}&\mathbf{a^{2}}&\mathbf{a^{3}}&\mathbf{a^{4}}&\mathbf{a^{5}}\\ \hline
\mathbf{       1}&   1&    &    &    &  -5&    &  -8&    &  -3&    \\
\mathbf{ z^{  2}}&    &   5&    &  11&    &  11&    &   9&    &   4\\
\mathbf{ z^{  4}}&  -8&    & -10&    &  15&    &  32&    &  15&    \\
\mathbf{ z^{  6}}&    & -18&    & -30&    & -22&    & -27&    & -17\\
\mathbf{ z^{  8}}&  14&    &  18&    & -15&    & -51&    & -32&    \\
\mathbf{ z^{ 10}}&    &  20&    &  27&    &  13&    &  26&    &  20\\
\mathbf{ z^{ 12}}&  -7&    &  -8&    &   7&    &  35&    &  27&    \\
\mathbf{ z^{ 14}}&    &  -8&    &  -9&    &  -2&    &  -9&    &  -8\\
\mathbf{ z^{ 16}}&   1&    &   1&    &  -1&    & -10&    &  -9&    \\
\mathbf{ z^{ 18}}&    &   1&    &   1&    &    &    &   1&    &   1\\
\mathbf{ z^{ 20}}&    &    &    &    &    &    &   1&    &   1&    \\
\hline
\end{array}\egroup$
}

\expandafter\def\csname HVolume_K14n_22589\endcsname{%
$8.878159662$
}

\expandafter\def\csname Alexander_K15n_51748\endcsname{%
\hbox{$1$}\allowbreak\thinspace
}

\expandafter\def\csname Jones_K15n_51748\endcsname{%
\hbox{$t^{-7}$}\allowbreak\thinspace\hbox{$-\,t^{-6}$}\allowbreak\thinspace\hbox{$-\,t^{-3}$}\allowbreak\thinspace\hbox{$+\,t^{-2}$}\allowbreak\thinspace\hbox{$-\,t^{-1}$}\allowbreak\thinspace\hbox{$+\,2$}\allowbreak\thinspace\hbox{$+\,t^{3}$}\allowbreak\thinspace\hbox{$-\,t^{4}$}\allowbreak\thinspace
}

\expandafter\def\csname HOMFLYPT_K15n_51748\endcsname{%
$\bgroup\arraycolsep=3pt\begin{array}{|r|rrrrr|}\hline
\vrule width 0pt depth 0pt height 10pt
&\mathbf{l^{-6}}&\mathbf{l^{-4}}&\mathbf{l^{-2}}&\mathbf{1}&\mathbf{l^{2}}\\ \hline
\mathbf{       1}&   1&    &  -6&   9&  -3\\
\mathbf{ m^{  2}}&   1&    & -11&  14&  -4\\
\mathbf{ m^{  4}}&    &    &  -6&   7&  -1\\
\mathbf{ m^{  6}}&    &    &  -1&   1&    \\
\hline
\end{array}\egroup$
}

\expandafter\def\csname Kauffman_K15n_51748\endcsname{%
$\bgroup\arraycolsep=1pt\begin{array}{|@{\hspace{3pt}}r@{\hspace{3pt}}|@{\hspace{3pt}}rrrrrrrrrr@{\hspace{3pt}}|}\hline
\vrule width 0pt depth 0pt height 10pt
&\mathbf{a^{-3}}&\mathbf{a^{-2}}&\mathbf{a^{-1}}&\mathbf{1}&\mathbf{a}&\mathbf{a^{2}}&\mathbf{a^{3}}&\mathbf{a^{4}}&\mathbf{a^{5}}&\mathbf{a^{6}}\\ \hline
\mathbf{       1}&    &   3&    &   9&    &   6&    &    &    &  -1\\
\mathbf{ z^{  2}}&  -3&    &  -9&    & -17&    & -19&    &  -8&    \\
\mathbf{ z^{  4}}&    & -10&    & -21&    &  -4&    &  19&    &  12\\
\mathbf{ z^{  6}}&   9&    &  17&    &  37&    &  62&    &  33&    \\
\mathbf{ z^{  8}}&    &  14&    &  21&    &  -7&    & -45&    & -31\\
\mathbf{ z^{ 10}}&  -6&    &  -8&    & -28&    & -78&    & -52&    \\
\mathbf{ z^{ 12}}&    &  -7&    &  -8&    &   6&    &  34&    &  27\\
\mathbf{ z^{ 14}}&   1&    &   1&    &   9&    &  44&    &  35&    \\
\mathbf{ z^{ 16}}&    &   1&    &   1&    &  -1&    & -10&    &  -9\\
\mathbf{ z^{ 18}}&    &    &    &    &  -1&    & -11&    & -10&    \\
\mathbf{ z^{ 20}}&    &    &    &    &    &    &    &   1&    &   1\\
\mathbf{ z^{ 22}}&    &    &    &    &    &    &   1&    &   1&    \\
\hline
\end{array}\egroup$
}

\expandafter\def\csname HVolume_K15n_51748\endcsname{%
$8.925447697$
}

\expandafter\def\csname Alexander_K15n_57436\endcsname{%
\hbox{$-\,t^{-2}$}\allowbreak\thinspace\hbox{$+\,3$}\allowbreak\thinspace\hbox{$-\,t^{2}$}\allowbreak\thinspace
}

\expandafter\def\csname Jones_K15n_57436\endcsname{%
\hbox{$t^{-7}$}\allowbreak\thinspace\hbox{$-\,t^{-6}$}\allowbreak\thinspace\hbox{$+\,t^{-4}$}\allowbreak\thinspace\hbox{$-\,2t^{-3}$}\allowbreak\thinspace\hbox{$+\,2t^{-2}$}\allowbreak\thinspace\hbox{$-\,2t^{-1}$}\allowbreak\thinspace\hbox{$+\,2$}\allowbreak\thinspace\hbox{$-\,t^{2}$}\allowbreak\thinspace\hbox{$+\,2t^{3}$}\allowbreak\thinspace\hbox{$-\,2t^{4}$}\allowbreak\thinspace\hbox{$+\,t^{5}$}\allowbreak\thinspace
}

\expandafter\def\csname HOMFLYPT_K15n_57436\endcsname{%
$\bgroup\arraycolsep=3pt\begin{array}{|r|rrrrrr|}\hline
\vrule width 0pt depth 0pt height 10pt
&\mathbf{l^{-6}}&\mathbf{l^{-4}}&\mathbf{l^{-2}}&\mathbf{1}&\mathbf{l^{2}}&\mathbf{l^{4}}\\ \hline
\mathbf{       1}&   2&  -4&   3&    &  -1&   1\\
\mathbf{ m^{  2}}&   1&  -5&   3&  -1&  -3&   1\\
\mathbf{ m^{  4}}&    &  -1&   1&    &  -1&    \\
\hline
\end{array}\egroup$
}

\expandafter\def\csname Kauffman_K15n_57436\endcsname{%
$\bgroup\arraycolsep=1pt\begin{array}{|@{\hspace{3pt}}r@{\hspace{3pt}}|@{\hspace{3pt}}rrrrrrrrrrr@{\hspace{3pt}}|}\hline
\vrule width 0pt depth 0pt height 10pt
&\mathbf{a^{-4}}&\mathbf{a^{-3}}&\mathbf{a^{-2}}&\mathbf{a^{-1}}&\mathbf{1}&\mathbf{a}&\mathbf{a^{2}}&\mathbf{a^{3}}&\mathbf{a^{4}}&\mathbf{a^{5}}&\mathbf{a^{6}}\\ \hline
\mathbf{       1}&   1&    &   1&    &    &    &  -3&    &  -4&    &  -2\\
\mathbf{ z^{  2}}&    &   2&    &   2&    &  -4&    & -10&    &  -6&    \\
\mathbf{ z^{  4}}&  -5&    &  -5&    &   8&    &  23&    &  30&    &  15\\
\mathbf{ z^{  6}}&    & -14&    & -18&    &  13&    &  44&    &  27&    \\
\mathbf{ z^{  8}}&  10&    &   4&    & -19&    & -40&    & -59&    & -32\\
\mathbf{ z^{ 10}}&    &  24&    &  31&    & -15&    & -69&    & -47&    \\
\mathbf{ z^{ 12}}&  -6&    &   7&    &  20&    &  21&    &  41&    &  27\\
\mathbf{ z^{ 14}}&    & -13&    & -15&    &   7&    &  43&    &  34&    \\
\mathbf{ z^{ 16}}&   1&    &  -6&    &  -8&    &  -3&    & -11&    &  -9\\
\mathbf{ z^{ 18}}&    &   2&    &   2&    &  -1&    & -11&    & -10&    \\
\mathbf{ z^{ 20}}&    &    &   1&    &   1&    &    &    &   1&    &   1\\
\mathbf{ z^{ 22}}&    &    &    &    &    &    &    &   1&    &   1&    \\
\hline
\end{array}\egroup$
}

\expandafter\def\csname HVolume_K15n_57436\endcsname{%
$12.529792456$
}

\expandafter\def\csname Alexander_K15n_57606\endcsname{%
\hbox{$-\,t^{-2}$}\allowbreak\thinspace\hbox{$+\,3$}\allowbreak\thinspace\hbox{$-\,t^{2}$}\allowbreak\thinspace
}

\expandafter\def\csname Jones_K15n_57606\endcsname{%
\hbox{$t^{-7}$}\allowbreak\thinspace\hbox{$-\,t^{-6}$}\allowbreak\thinspace\hbox{$+\,t^{-4}$}\allowbreak\thinspace\hbox{$-\,2t^{-3}$}\allowbreak\thinspace\hbox{$+\,2t^{-2}$}\allowbreak\thinspace\hbox{$-\,2t^{-1}$}\allowbreak\thinspace\hbox{$+\,2$}\allowbreak\thinspace\hbox{$-\,t^{2}$}\allowbreak\thinspace\hbox{$+\,2t^{3}$}\allowbreak\thinspace\hbox{$-\,2t^{4}$}\allowbreak\thinspace\hbox{$+\,t^{5}$}\allowbreak\thinspace
}

\expandafter\def\csname HOMFLYPT_K15n_57606\endcsname{%
$\bgroup\arraycolsep=3pt\begin{array}{|r|rrrrrr|}\hline
\vrule width 0pt depth 0pt height 10pt
&\mathbf{l^{-6}}&\mathbf{l^{-4}}&\mathbf{l^{-2}}&\mathbf{1}&\mathbf{l^{2}}&\mathbf{l^{4}}\\ \hline
\mathbf{       1}&   2&  -4&   3&    &  -1&   1\\
\mathbf{ m^{  2}}&   1&  -5&   3&  -1&  -3&   1\\
\mathbf{ m^{  4}}&    &  -1&   1&    &  -1&    \\
\hline
\end{array}\egroup$
}

\expandafter\def\csname Kauffman_K15n_57606\endcsname{%
$\bgroup\arraycolsep=1pt\begin{array}{|@{\hspace{3pt}}r@{\hspace{3pt}}|@{\hspace{3pt}}rrrrrrrrrrr@{\hspace{3pt}}|}\hline
\vrule width 0pt depth 0pt height 10pt
&\mathbf{a^{-4}}&\mathbf{a^{-3}}&\mathbf{a^{-2}}&\mathbf{a^{-1}}&\mathbf{1}&\mathbf{a}&\mathbf{a^{2}}&\mathbf{a^{3}}&\mathbf{a^{4}}&\mathbf{a^{5}}&\mathbf{a^{6}}\\ \hline
\mathbf{       1}&   1&    &   1&    &    &    &  -3&    &  -4&    &  -2\\
\mathbf{ z^{  2}}&    &   2&    &   2&    &  -4&    & -10&    &  -6&    \\
\mathbf{ z^{  4}}&  -5&    &  -5&    &   8&    &  23&    &  30&    &  15\\
\mathbf{ z^{  6}}&    & -14&    & -18&    &  13&    &  44&    &  27&    \\
\mathbf{ z^{  8}}&  10&    &   4&    & -19&    & -40&    & -59&    & -32\\
\mathbf{ z^{ 10}}&    &  24&    &  31&    & -15&    & -69&    & -47&    \\
\mathbf{ z^{ 12}}&  -6&    &   7&    &  20&    &  21&    &  41&    &  27\\
\mathbf{ z^{ 14}}&    & -13&    & -15&    &   7&    &  43&    &  34&    \\
\mathbf{ z^{ 16}}&   1&    &  -6&    &  -8&    &  -3&    & -11&    &  -9\\
\mathbf{ z^{ 18}}&    &   2&    &   2&    &  -1&    & -11&    & -10&    \\
\mathbf{ z^{ 20}}&    &    &   1&    &   1&    &    &    &   1&    &   1\\
\mathbf{ z^{ 22}}&    &    &    &    &    &    &    &   1&    &   1&    \\
\hline
\end{array}\egroup$
}

\expandafter\def\csname HVolume_K15n_57606\endcsname{%
$12.529792456$
}

\expandafter\def\csname Alexander_K15n_115375\endcsname{%
\hbox{$1$}\allowbreak\thinspace
}

\expandafter\def\csname Jones_K15n_115375\endcsname{%
\hbox{$t^{-7}$}\allowbreak\thinspace\hbox{$-\,t^{-6}$}\allowbreak\thinspace\hbox{$-\,t^{-3}$}\allowbreak\thinspace\hbox{$+\,t^{-2}$}\allowbreak\thinspace\hbox{$-\,t^{-1}$}\allowbreak\thinspace\hbox{$+\,2$}\allowbreak\thinspace\hbox{$+\,t^{3}$}\allowbreak\thinspace\hbox{$-\,t^{4}$}\allowbreak\thinspace
}

\expandafter\def\csname HOMFLYPT_K15n_115375\endcsname{%
$\bgroup\arraycolsep=3pt\begin{array}{|r|rrrrr|}\hline
\vrule width 0pt depth 0pt height 10pt
&\mathbf{l^{-6}}&\mathbf{l^{-4}}&\mathbf{l^{-2}}&\mathbf{1}&\mathbf{l^{2}}\\ \hline
\mathbf{       1}&   1&    &  -6&   9&  -3\\
\mathbf{ m^{  2}}&   1&    & -11&  14&  -4\\
\mathbf{ m^{  4}}&    &    &  -6&   7&  -1\\
\mathbf{ m^{  6}}&    &    &  -1&   1&    \\
\hline
\end{array}\egroup$
}

\expandafter\def\csname Kauffman_K15n_115375\endcsname{%
$\bgroup\arraycolsep=1pt\begin{array}{|@{\hspace{3pt}}r@{\hspace{3pt}}|@{\hspace{3pt}}rrrrrrrrrr@{\hspace{3pt}}|}\hline
\vrule width 0pt depth 0pt height 10pt
&\mathbf{a^{-3}}&\mathbf{a^{-2}}&\mathbf{a^{-1}}&\mathbf{1}&\mathbf{a}&\mathbf{a^{2}}&\mathbf{a^{3}}&\mathbf{a^{4}}&\mathbf{a^{5}}&\mathbf{a^{6}}\\ \hline
\mathbf{       1}&    &   3&    &   9&    &   6&    &    &    &  -1\\
\mathbf{ z^{  2}}&  -3&    &  -9&    & -17&    & -19&    &  -8&    \\
\mathbf{ z^{  4}}&    & -10&    & -21&    &  -4&    &  19&    &  12\\
\mathbf{ z^{  6}}&   9&    &  17&    &  37&    &  62&    &  33&    \\
\mathbf{ z^{  8}}&    &  14&    &  21&    &  -7&    & -45&    & -31\\
\mathbf{ z^{ 10}}&  -6&    &  -8&    & -28&    & -78&    & -52&    \\
\mathbf{ z^{ 12}}&    &  -7&    &  -8&    &   6&    &  34&    &  27\\
\mathbf{ z^{ 14}}&   1&    &   1&    &   9&    &  44&    &  35&    \\
\mathbf{ z^{ 16}}&    &   1&    &   1&    &  -1&    & -10&    &  -9\\
\mathbf{ z^{ 18}}&    &    &    &    &  -1&    & -11&    & -10&    \\
\mathbf{ z^{ 20}}&    &    &    &    &    &    &    &   1&    &   1\\
\mathbf{ z^{ 22}}&    &    &    &    &    &    &   1&    &   1&    \\
\hline
\end{array}\egroup$
}

\expandafter\def\csname HVolume_K15n_115375\endcsname{%
$8.925447697$
}

\expandafter\def\csname Alexander_K15n_133697\endcsname{%
\hbox{$t^{-3}$}\allowbreak\thinspace\hbox{$-\,t^{-2}$}\allowbreak\thinspace\hbox{$-\,t^{-1}$}\allowbreak\thinspace\hbox{$+\,3$}\allowbreak\thinspace\hbox{$-\,t$}\allowbreak\thinspace\hbox{$-\,t^{2}$}\allowbreak\thinspace\hbox{$+\,t^{3}$}\allowbreak\thinspace
}

\expandafter\def\csname Jones_K15n_133697\endcsname{%
\hbox{$-\,t^{-6}$}\allowbreak\thinspace\hbox{$+\,2t^{-5}$}\allowbreak\thinspace\hbox{$-\,2t^{-4}$}\allowbreak\thinspace\hbox{$+\,t^{-3}$}\allowbreak\thinspace\hbox{$-\,t^{-1}$}\allowbreak\thinspace\hbox{$+\,3$}\allowbreak\thinspace\hbox{$-\,2t$}\allowbreak\thinspace\hbox{$+\,2t^{2}$}\allowbreak\thinspace\hbox{$-\,t^{3}$}\allowbreak\thinspace\hbox{$+\,t^{5}$}\allowbreak\thinspace\hbox{$-\,t^{6}$}\allowbreak\thinspace
}

\expandafter\def\csname HOMFLYPT_K15n_133697\endcsname{%
$\bgroup\arraycolsep=3pt\begin{array}{|r|rrrrr|}\hline
\vrule width 0pt depth 0pt height 10pt
&\mathbf{l^{-4}}&\mathbf{l^{-2}}&\mathbf{1}&\mathbf{l^{2}}&\mathbf{l^{4}}\\ \hline
\mathbf{       1}&  -2&   5&  -6&   7&  -3\\
\mathbf{ m^{  2}}&  -3&  10& -12&  13&  -4\\
\mathbf{ m^{  4}}&  -1&   6&  -6&   7&  -1\\
\mathbf{ m^{  6}}&    &   1&  -1&   1&    \\
\hline
\end{array}\egroup$
}

\expandafter\def\csname Kauffman_K15n_133697\endcsname{%
$\bgroup\arraycolsep=1pt\begin{array}{|@{\hspace{3pt}}r@{\hspace{3pt}}|@{\hspace{3pt}}rrrrrrrrrrr@{\hspace{3pt}}|}\hline
\vrule width 0pt depth 0pt height 10pt
&\mathbf{a^{-5}}&\mathbf{a^{-4}}&\mathbf{a^{-3}}&\mathbf{a^{-2}}&\mathbf{a^{-1}}&\mathbf{1}&\mathbf{a}&\mathbf{a^{2}}&\mathbf{a^{3}}&\mathbf{a^{4}}&\mathbf{a^{5}}\\ \hline
\mathbf{       1}&    &  -3&    &  -7&    &  -6&    &  -5&    &  -2&    \\
\mathbf{ z^{  2}}&   6&    &  11&    &   6&    &   2&    &   4&    &   3\\
\mathbf{ z^{  4}}&    &  12&    &  28&    &  30&    &  23&    &   9&    \\
\mathbf{ z^{  6}}& -18&    & -35&    & -10&    &   6&    & -12&    & -11\\
\mathbf{ z^{  8}}&    & -28&    & -48&    & -39&    & -48&    & -29&    \\
\mathbf{ z^{ 10}}&  20&    &  32&    &  -1&    & -26&    &   2&    &  15\\
\mathbf{ z^{ 12}}&    &  26&    &  35&    &  20&    &  46&    &  35&    \\
\mathbf{ z^{ 14}}&  -8&    & -10&    &   5&    &  26&    &  12&    &  -7\\
\mathbf{ z^{ 16}}&    &  -9&    & -10&    &  -3&    & -17&    & -15&    \\
\mathbf{ z^{ 18}}&   1&    &   1&    &  -1&    &  -9&    &  -7&    &   1\\
\mathbf{ z^{ 20}}&    &   1&    &   1&    &    &    &   2&    &   2&    \\
\mathbf{ z^{ 22}}&    &    &    &    &    &    &   1&    &   1&    &    \\
\hline
\end{array}\egroup$
}

\expandafter\def\csname HVolume_K15n_133697\endcsname{%
$12.569864535$
}

\expandafter\def\csname Alexander_K15n_135711\endcsname{%
\hbox{$t^{-3}$}\allowbreak\thinspace\hbox{$-\,t^{-2}$}\allowbreak\thinspace\hbox{$-\,t^{-1}$}\allowbreak\thinspace\hbox{$+\,3$}\allowbreak\thinspace\hbox{$-\,t$}\allowbreak\thinspace\hbox{$-\,t^{2}$}\allowbreak\thinspace\hbox{$+\,t^{3}$}\allowbreak\thinspace
}

\expandafter\def\csname Jones_K15n_135711\endcsname{%
\hbox{$-\,t^{-6}$}\allowbreak\thinspace\hbox{$+\,2t^{-5}$}\allowbreak\thinspace\hbox{$-\,2t^{-4}$}\allowbreak\thinspace\hbox{$+\,t^{-3}$}\allowbreak\thinspace\hbox{$-\,t^{-1}$}\allowbreak\thinspace\hbox{$+\,3$}\allowbreak\thinspace\hbox{$-\,2t$}\allowbreak\thinspace\hbox{$+\,2t^{2}$}\allowbreak\thinspace\hbox{$-\,t^{3}$}\allowbreak\thinspace\hbox{$+\,t^{5}$}\allowbreak\thinspace\hbox{$-\,t^{6}$}\allowbreak\thinspace
}

\expandafter\def\csname HOMFLYPT_K15n_135711\endcsname{%
$\bgroup\arraycolsep=3pt\begin{array}{|r|rrrrr|}\hline
\vrule width 0pt depth 0pt height 10pt
&\mathbf{l^{-4}}&\mathbf{l^{-2}}&\mathbf{1}&\mathbf{l^{2}}&\mathbf{l^{4}}\\ \hline
\mathbf{       1}&  -2&   5&  -6&   7&  -3\\
\mathbf{ m^{  2}}&  -3&  10& -12&  13&  -4\\
\mathbf{ m^{  4}}&  -1&   6&  -6&   7&  -1\\
\mathbf{ m^{  6}}&    &   1&  -1&   1&    \\
\hline
\end{array}\egroup$
}

\expandafter\def\csname Kauffman_K15n_135711\endcsname{%
$\bgroup\arraycolsep=1pt\begin{array}{|@{\hspace{3pt}}r@{\hspace{3pt}}|@{\hspace{3pt}}rrrrrrrrrrr@{\hspace{3pt}}|}\hline
\vrule width 0pt depth 0pt height 10pt
&\mathbf{a^{-5}}&\mathbf{a^{-4}}&\mathbf{a^{-3}}&\mathbf{a^{-2}}&\mathbf{a^{-1}}&\mathbf{1}&\mathbf{a}&\mathbf{a^{2}}&\mathbf{a^{3}}&\mathbf{a^{4}}&\mathbf{a^{5}}\\ \hline
\mathbf{       1}&    &  -3&    &  -7&    &  -6&    &  -5&    &  -2&    \\
\mathbf{ z^{  2}}&   6&    &  11&    &   6&    &   2&    &   4&    &   3\\
\mathbf{ z^{  4}}&    &  12&    &  28&    &  30&    &  23&    &   9&    \\
\mathbf{ z^{  6}}& -18&    & -35&    & -10&    &   6&    & -12&    & -11\\
\mathbf{ z^{  8}}&    & -28&    & -48&    & -39&    & -48&    & -29&    \\
\mathbf{ z^{ 10}}&  20&    &  32&    &  -1&    & -26&    &   2&    &  15\\
\mathbf{ z^{ 12}}&    &  26&    &  35&    &  20&    &  46&    &  35&    \\
\mathbf{ z^{ 14}}&  -8&    & -10&    &   5&    &  26&    &  12&    &  -7\\
\mathbf{ z^{ 16}}&    &  -9&    & -10&    &  -3&    & -17&    & -15&    \\
\mathbf{ z^{ 18}}&   1&    &   1&    &  -1&    &  -9&    &  -7&    &   1\\
\mathbf{ z^{ 20}}&    &   1&    &   1&    &    &    &   2&    &   2&    \\
\mathbf{ z^{ 22}}&    &    &    &    &    &    &   1&    &   1&    &    \\
\hline
\end{array}\egroup$
}

\expandafter\def\csname HVolume_K15n_135711\endcsname{%
$12.569864535$
}

\expandafter\def\csname Alexander_K15n_148673\endcsname{%
\hbox{$t^{-4}$}\allowbreak\thinspace\hbox{$-\,2t^{-3}$}\allowbreak\thinspace\hbox{$+\,t^{-2}$}\allowbreak\thinspace\hbox{$+\,3t^{-1}$}\allowbreak\thinspace\hbox{$-\,5$}\allowbreak\thinspace\hbox{$+\,3t$}\allowbreak\thinspace\hbox{$+\,t^{2}$}\allowbreak\thinspace\hbox{$-\,2t^{3}$}\allowbreak\thinspace\hbox{$+\,t^{4}$}\allowbreak\thinspace
}

\expandafter\def\csname Jones_K15n_148673\endcsname{%
\hbox{$-\,t^{-3}$}\allowbreak\thinspace\hbox{$+\,2t^{-2}$}\allowbreak\thinspace\hbox{$-\,2t^{-1}$}\allowbreak\thinspace\hbox{$+\,1$}\allowbreak\thinspace\hbox{$+\,t$}\allowbreak\thinspace\hbox{$-\,t^{2}$}\allowbreak\thinspace\hbox{$+\,3t^{3}$}\allowbreak\thinspace\hbox{$-\,3t^{4}$}\allowbreak\thinspace\hbox{$+\,2t^{5}$}\allowbreak\thinspace\hbox{$-\,t^{6}$}\allowbreak\thinspace\hbox{$+\,t^{8}$}\allowbreak\thinspace\hbox{$-\,t^{9}$}\allowbreak\thinspace
}

\expandafter\def\csname HOMFLYPT_K15n_148673\endcsname{%
$\bgroup\arraycolsep=3pt\begin{array}{|r|rrrrr|}\hline
\vrule width 0pt depth 0pt height 10pt
&\mathbf{1}&\mathbf{l^{2}}&\mathbf{l^{4}}&\mathbf{l^{6}}&\mathbf{l^{8}}\\ \hline
\mathbf{       1}&  -2&   5&  -2&   1&  -1\\
\mathbf{ m^{  2}}&  -6&  13&  -2&   1&  -1\\
\mathbf{ m^{  4}}&  -5&  15&  -1&    &    \\
\mathbf{ m^{  6}}&  -1&   7&    &    &    \\
\mathbf{ m^{  8}}&    &   1&    &    &    \\
\hline
\end{array}\egroup$
}

\expandafter\def\csname Kauffman_K15n_148673\endcsname{%
$\bgroup\arraycolsep=1pt\begin{array}{|@{\hspace{3pt}}r@{\hspace{3pt}}|@{\hspace{3pt}}rrrrrrrrrrr@{\hspace{3pt}}|}\hline
\vrule width 0pt depth 0pt height 10pt
&\mathbf{a^{-9}}&\mathbf{a^{-8}}&\mathbf{a^{-7}}&\mathbf{a^{-6}}&\mathbf{a^{-5}}&\mathbf{a^{-4}}&\mathbf{a^{-3}}&\mathbf{a^{-2}}&\mathbf{a^{-1}}&\mathbf{1}&\mathbf{a}\\ \hline
\mathbf{       1}&    &  -1&    &  -1&    &  -2&    &  -5&    &  -2&    \\
\mathbf{ z^{  2}}&   5&    &   5&    &  -2&    &  -2&    &   2&    &   2\\
\mathbf{ z^{  4}}&    &  10&    &  12&    &  13&    &  22&    &  11&    \\
\mathbf{ z^{  6}}& -18&    & -22&    &  10&    &  16&    &  -9&    & -11\\
\mathbf{ z^{  8}}&    & -28&    & -32&    & -22&    & -48&    & -30&    \\
\mathbf{ z^{ 10}}&  20&    &  25&    & -13&    & -32&    &   1&    &  15\\
\mathbf{ z^{ 12}}&    &  26&    &  28&    &  13&    &  46&    &  35&    \\
\mathbf{ z^{ 14}}&  -8&    &  -9&    &   7&    &  27&    &  12&    &  -7\\
\mathbf{ z^{ 16}}&    &  -9&    &  -9&    &  -2&    & -17&    & -15&    \\
\mathbf{ z^{ 18}}&   1&    &   1&    &  -1&    &  -9&    &  -7&    &   1\\
\mathbf{ z^{ 20}}&    &   1&    &   1&    &    &    &   2&    &   2&    \\
\mathbf{ z^{ 22}}&    &    &    &    &    &    &   1&    &   1&    &    \\
\hline
\end{array}\egroup$
}

\expandafter\def\csname HVolume_K15n_148673\endcsname{%
$13.081220984$
}

\expandafter\def\csname Alexander_K15n_151500M\endcsname{%
\hbox{$t^{-4}$}\allowbreak\thinspace\hbox{$-\,2t^{-3}$}\allowbreak\thinspace\hbox{$+\,t^{-2}$}\allowbreak\thinspace\hbox{$+\,3t^{-1}$}\allowbreak\thinspace\hbox{$-\,5$}\allowbreak\thinspace\hbox{$+\,3t$}\allowbreak\thinspace\hbox{$+\,t^{2}$}\allowbreak\thinspace\hbox{$-\,2t^{3}$}\allowbreak\thinspace\hbox{$+\,t^{4}$}\allowbreak\thinspace
}

\expandafter\def\csname Jones_K15n_151500M\endcsname{%
\hbox{$-\,t^{-3}$}\allowbreak\thinspace\hbox{$+\,2t^{-2}$}\allowbreak\thinspace\hbox{$-\,2t^{-1}$}\allowbreak\thinspace\hbox{$+\,1$}\allowbreak\thinspace\hbox{$+\,t$}\allowbreak\thinspace\hbox{$-\,t^{2}$}\allowbreak\thinspace\hbox{$+\,3t^{3}$}\allowbreak\thinspace\hbox{$-\,3t^{4}$}\allowbreak\thinspace\hbox{$+\,2t^{5}$}\allowbreak\thinspace\hbox{$-\,t^{6}$}\allowbreak\thinspace\hbox{$+\,t^{8}$}\allowbreak\thinspace\hbox{$-\,t^{9}$}\allowbreak\thinspace
}

\expandafter\def\csname HOMFLYPT_K15n_151500M\endcsname{%
$\bgroup\arraycolsep=3pt\begin{array}{|r|rrrrr|}\hline
\vrule width 0pt depth 0pt height 10pt
&\mathbf{1}&\mathbf{l^{2}}&\mathbf{l^{4}}&\mathbf{l^{6}}&\mathbf{l^{8}}\\ \hline
\mathbf{       1}&  -2&   5&  -2&   1&  -1\\
\mathbf{ m^{  2}}&  -6&  13&  -2&   1&  -1\\
\mathbf{ m^{  4}}&  -5&  15&  -1&    &    \\
\mathbf{ m^{  6}}&  -1&   7&    &    &    \\
\mathbf{ m^{  8}}&    &   1&    &    &    \\
\hline
\end{array}\egroup$
}

\expandafter\def\csname Kauffman_K15n_151500M\endcsname{%
$\bgroup\arraycolsep=1pt\begin{array}{|@{\hspace{3pt}}r@{\hspace{3pt}}|@{\hspace{3pt}}rrrrrrrrrrr@{\hspace{3pt}}|}\hline
\vrule width 0pt depth 0pt height 10pt
&\mathbf{a^{-9}}&\mathbf{a^{-8}}&\mathbf{a^{-7}}&\mathbf{a^{-6}}&\mathbf{a^{-5}}&\mathbf{a^{-4}}&\mathbf{a^{-3}}&\mathbf{a^{-2}}&\mathbf{a^{-1}}&\mathbf{1}&\mathbf{a}\\ \hline
\mathbf{       1}&    &  -1&    &  -1&    &  -2&    &  -5&    &  -2&    \\
\mathbf{ z^{  2}}&   5&    &   5&    &  -2&    &  -2&    &   2&    &   2\\
\mathbf{ z^{  4}}&    &  10&    &  12&    &  13&    &  22&    &  11&    \\
\mathbf{ z^{  6}}& -18&    & -22&    &  10&    &  16&    &  -9&    & -11\\
\mathbf{ z^{  8}}&    & -28&    & -32&    & -22&    & -48&    & -30&    \\
\mathbf{ z^{ 10}}&  20&    &  25&    & -13&    & -32&    &   1&    &  15\\
\mathbf{ z^{ 12}}&    &  26&    &  28&    &  13&    &  46&    &  35&    \\
\mathbf{ z^{ 14}}&  -8&    &  -9&    &   7&    &  27&    &  12&    &  -7\\
\mathbf{ z^{ 16}}&    &  -9&    &  -9&    &  -2&    & -17&    & -15&    \\
\mathbf{ z^{ 18}}&   1&    &   1&    &  -1&    &  -9&    &  -7&    &   1\\
\mathbf{ z^{ 20}}&    &   1&    &   1&    &    &    &   2&    &   2&    \\
\mathbf{ z^{ 22}}&    &    &    &    &    &    &   1&    &   1&    &    \\
\hline
\end{array}\egroup$
}

\expandafter\def\csname HVolume_K15n_151500M\endcsname{%
$13.081220984$
}

\expandafter\def\csname KhovanovQ_K14n_22185\endcsname{%
$1_{\underline{13}}^{\underline{7}}\, 1_{\underline{9}}^{\underline{6}}\, 1_{\underline{7}}^{\underline{4}}\, 1_{\underline{7}}^{\underline{3}}\, 1_{\underline{3}}^{\underline{3}}\, 1_{\underline{5}}^{\underline{2}}\, 1_{\underline{3}}^{\underline{2}}\, 1_{\underline{3}}^{\underline{1}}\, 1_{\underline{1}}^{\underline{1}}\, 1_{\underline{3}}^{0}\, 2_{\underline{1}}^{0}\, 2_{1}^{0}\, 2_{1}^{1}\, 1_{3}^{1}\, 1_{1}^{2}\, 1_{3}^{2}\, 1_{5}^{2}\, 1_{3}^{3}\, 1_{5}^{3}\, 1_{7}^{3}\, 1_{7}^{4}\, 1_{7}^{5}\, 1_{11}^{6}$
}

\expandafter\def\csname KhovanovQ_K14n_22589\endcsname{%
$1_{\underline{13}}^{\underline{7}}\, 1_{\underline{9}}^{\underline{6}}\, 1_{\underline{9}}^{\underline{5}}\, 1_{\underline{9}}^{\underline{4}}\, 1_{\underline{7}}^{\underline{4}}\, 1_{\underline{5}}^{\underline{4}}\, 1_{\underline{7}}^{\underline{3}}\, 1_{\underline{5}}^{\underline{3}}\, 1_{\underline{3}}^{\underline{3}}\, 1_{\underline{5}}^{\underline{2}}\, 2_{\underline{3}}^{\underline{2}}\, 1_{\underline{3}}^{\underline{1}}\, 1_{\underline{1}}^{\underline{1}}\, 1_{1}^{\underline{1}}\, 2_{\underline{1}}^{0}\, 2_{1}^{0}\, 1_{1}^{1}\, 1_{3}^{1}\, 1_{1}^{2}\, 1_{5}^{2}\, 1_{5}^{3}\, 1_{7}^{5}\, 1_{11}^{6}$
}

\expandafter\def\csname KhovanovQ_K15n_51748\endcsname{%
$1_{\underline{15}}^{\underline{8}}\, 1_{\underline{11}}^{\underline{7}}\, 1_{\underline{11}}^{\underline{6}}\, 1_{\underline{11}}^{\underline{5}}\, 1_{\underline{9}}^{\underline{5}}\, 1_{\underline{7}}^{\underline{5}}\, 1_{\underline{9}}^{\underline{4}}\, 1_{\underline{7}}^{\underline{4}}\, 1_{\underline{5}}^{\underline{4}}\, 1_{\underline{7}}^{\underline{3}}\, 2_{\underline{5}}^{\underline{3}}\, 1_{\underline{5}}^{\underline{2}}\, 1_{\underline{3}}^{\underline{2}}\, 1_{\underline{1}}^{\underline{2}}\, 1_{\underline{3}}^{\underline{1}}\, 1_{\underline{1}}^{\underline{1}}\, 2_{\underline{1}}^{0}\, 2_{1}^{0}\, 1_{\underline{1}}^{1}\, 1_{3}^{1}\, 1_{3}^{2}\, 1_{5}^{4}\, 1_{9}^{5}$
}

\expandafter\def\csname KhovanovQ_K15n_57436\endcsname{%
$1_{\underline{15}}^{\underline{8}}\, 1_{\underline{11}}^{\underline{7}}\, 1_{\underline{9}}^{\underline{5}}\, 2_{\underline{9}}^{\underline{4}}\, 1_{\underline{5}}^{\underline{4}}\, 1_{\underline{7}}^{\underline{3}}\, 2_{\underline{5}}^{\underline{3}}\, 2_{\underline{5}}^{\underline{2}}\, 1_{\underline{3}}^{\underline{2}}\, 1_{\underline{5}}^{\underline{1}}\, 2_{\underline{3}}^{\underline{1}}\, 2_{\underline{1}}^{\underline{1}}\, 1_{\underline{3}}^{0}\, 3_{\underline{1}}^{0}\, 3_{1}^{0}\, 1_{\underline{1}}^{1}\, 3_{1}^{1}\, 1_{3}^{1}\, 2_{1}^{2}\, 1_{3}^{2}\, 2_{5}^{2}\, 1_{3}^{3}\, 2_{5}^{3}\, 1_{5}^{4}\, 1_{7}^{4}\, 1_{7}^{5}\, 1_{9}^{5}\, 1_{11}^{6}$
}

\expandafter\def\csname KhovanovQ_K15n_57606\endcsname{%
$1_{\underline{15}}^{\underline{8}}\, 1_{\underline{11}}^{\underline{7}}\, 1_{\underline{11}}^{\underline{6}}\, 1_{\underline{11}}^{\underline{5}}\, 1_{\underline{9}}^{\underline{5}}\, 1_{\underline{7}}^{\underline{5}}\, 2_{\underline{9}}^{\underline{4}}\, 1_{\underline{7}}^{\underline{4}}\, 1_{\underline{5}}^{\underline{4}}\, 1_{\underline{7}}^{\underline{3}}\, 3_{\underline{5}}^{\underline{3}}\, 2_{\underline{5}}^{\underline{2}}\, 1_{\underline{3}}^{\underline{2}}\, 1_{\underline{1}}^{\underline{2}}\, 2_{\underline{3}}^{\underline{1}}\, 2_{\underline{1}}^{\underline{1}}\, 1_{\underline{3}}^{0}\, 2_{\underline{1}}^{0}\, 3_{1}^{0}\, 1_{\underline{1}}^{1}\, 2_{1}^{1}\, 1_{3}^{1}\, 1_{1}^{2}\, 1_{3}^{2}\, 1_{5}^{2}\, 1_{3}^{3}\, 1_{5}^{3}\, 1_{5}^{4}\, 1_{7}^{4}\, 1_{7}^{5}\, 1_{9}^{5}\, 1_{11}^{6}$
}

\expandafter\def\csname KhovanovQ_K15n_115375\endcsname{%
$1_{\underline{15}}^{\underline{8}}\, 1_{\underline{11}}^{\underline{7}}\, 1_{\underline{9}}^{\underline{5}}\, 1_{\underline{9}}^{\underline{4}}\, 1_{\underline{5}}^{\underline{4}}\, 1_{\underline{7}}^{\underline{3}}\, 1_{\underline{5}}^{\underline{3}}\, 1_{\underline{5}}^{\underline{2}}\, 1_{\underline{3}}^{\underline{2}}\, 1_{\underline{5}}^{\underline{1}}\, 1_{\underline{3}}^{\underline{1}}\, 1_{\underline{1}}^{\underline{1}}\, 3_{\underline{1}}^{0}\, 2_{1}^{0}\, 1_{\underline{1}}^{1}\, 1_{1}^{1}\, 1_{3}^{1}\, 1_{1}^{2}\, 1_{3}^{2}\, 1_{5}^{2}\, 1_{5}^{3}\, 1_{5}^{4}\, 1_{9}^{5}$
}

\expandafter\def\csname KhovanovQ_K15n_133697\endcsname{%
$1_{\underline{13}}^{\underline{7}}\, 1_{\underline{11}}^{\underline{6}}\, 1_{\underline{9}}^{\underline{6}}\, 1_{\underline{9}}^{\underline{5}}\, 1_{\underline{7}}^{\underline{5}}\, 1_{\underline{7}}^{\underline{4}}\, 1_{\underline{5}}^{\underline{4}}\, 1_{\underline{7}}^{\underline{3}}\, 1_{\underline{5}}^{\underline{3}}\, 1_{\underline{3}}^{\underline{3}}\, 1_{\underline{5}}^{\underline{2}}\, 2_{\underline{3}}^{\underline{2}}\, 1_{\underline{1}}^{\underline{2}}\, 2_{\underline{3}}^{\underline{1}}\, 1_{\underline{1}}^{\underline{1}}\, 1_{1}^{\underline{1}}\, 3_{\underline{1}}^{0}\, 3_{1}^{0}\, 1_{\underline{1}}^{1}\, 1_{1}^{1}\, 2_{3}^{1}\, 3_{3}^{2}\, 1_{5}^{2}\, 1_{3}^{3}\, 1_{5}^{3}\, 2_{7}^{3}\, 1_{5}^{4}\, 1_{7}^{4}\, 1_{9}^{4}\, 1_{9}^{5}\, 1_{9}^{6}\, 1_{13}^{7}$
}

\expandafter\def\csname KhovanovQ_K15n_135711\endcsname{%
$1_{\underline{13}}^{\underline{7}}\, 1_{\underline{11}}^{\underline{6}}\, 1_{\underline{9}}^{\underline{6}}\, 1_{\underline{9}}^{\underline{5}}\, 1_{\underline{7}}^{\underline{5}}\, 2_{\underline{7}}^{\underline{4}}\, 1_{\underline{5}}^{\underline{4}}\, 2_{\underline{7}}^{\underline{3}}\, 1_{\underline{5}}^{\underline{3}}\, 2_{\underline{3}}^{\underline{3}}\, 1_{\underline{5}}^{\underline{2}}\, 3_{\underline{3}}^{\underline{2}}\, 1_{\underline{1}}^{\underline{2}}\, 2_{\underline{3}}^{\underline{1}}\, 2_{\underline{1}}^{\underline{1}}\, 1_{1}^{\underline{1}}\, 3_{\underline{1}}^{0}\, 3_{1}^{0}\, 1_{3}^{0}\, 1_{1}^{1}\, 2_{3}^{1}\, 2_{3}^{2}\, 1_{5}^{2}\, 1_{3}^{3}\, 2_{7}^{3}\, 1_{7}^{4}\, 1_{9}^{6}\, 1_{13}^{7}$
}

\expandafter\def\csname KhovanovQ_K15n_148673\endcsname{%
$1_{\underline{7}}^{\underline{5}}\, 1_{\underline{5}}^{\underline{4}}\, 1_{\underline{3}}^{\underline{4}}\, 1_{\underline{3}}^{\underline{3}}\, 1_{\underline{1}}^{\underline{3}}\, 1_{\underline{1}}^{\underline{2}}\, 1_{1}^{\underline{2}}\, 1_{\underline{1}}^{\underline{1}}\, 1_{1}^{\underline{1}}\, 1_{3}^{\underline{1}}\, 2_{1}^{0}\, 3_{3}^{0}\, 1_{5}^{0}\, 2_{3}^{1}\, 1_{5}^{1}\, 1_{7}^{1}\, 3_{5}^{2}\, 2_{7}^{2}\, 1_{5}^{3}\, 1_{7}^{3}\, 3_{9}^{3}\, 3_{9}^{4}\, 1_{11}^{4}\, 1_{9}^{5}\, 1_{11}^{5}\, 2_{13}^{5}\, 1_{11}^{6}\, 1_{13}^{6}\, 1_{15}^{6}\, 1_{15}^{7}\, 1_{15}^{8}\, 1_{19}^{9}$
}

\expandafter\def\csname KhovanovQ_K15n_151500M\endcsname{%
$1_{\underline{7}}^{\underline{5}}\, 1_{\underline{5}}^{\underline{4}}\, 1_{\underline{3}}^{\underline{4}}\, 1_{\underline{3}}^{\underline{3}}\, 1_{\underline{1}}^{\underline{3}}\, 2_{\underline{1}}^{\underline{2}}\, 1_{1}^{\underline{2}}\, 2_{\underline{1}}^{\underline{1}}\, 1_{1}^{\underline{1}}\, 2_{3}^{\underline{1}}\, 2_{1}^{0}\, 4_{3}^{0}\, 1_{5}^{0}\, 2_{3}^{1}\, 2_{5}^{1}\, 1_{7}^{1}\, 3_{5}^{2}\, 2_{7}^{2}\, 1_{9}^{2}\, 1_{7}^{3}\, 3_{9}^{3}\, 2_{9}^{4}\, 1_{11}^{4}\, 1_{9}^{5}\, 2_{13}^{5}\, 1_{13}^{6}\, 1_{15}^{8}\, 1_{19}^{9}$
}

\expandafter\def\csname KhovanovT_K14n_22185\endcsname{%
$1_{\underline{11}}^{\underline{6}}\, 1_{\underline{9}}^{\underline{4}}\, 1_{\underline{7}}^{\underline{4}}\, 1_{\underline{7}}^{\underline{3}}\, 2_{\underline{5}}^{\underline{3}}\, 1_{\underline{5}}^{\underline{2}}\, 2_{\underline{3}}^{\underline{1}}\, 1_{\underline{1}}^{\underline{1}}\, 2_{\underline{1}}^{0}\, 1_{1}^{0}\, 1_{\underline{1}}^{1}\, 1_{1}^{1}\, 1_{3}^{2}\, 1_{3}^{3}\, 1_{5}^{3}\, 1_{5}^{4}\, 1_{9}^{6}$
}

\expandafter\def\csname KhovanovT_K14n_22589\endcsname{%
$1_{\underline{11}}^{\underline{6}}\, 1_{\underline{7}}^{\underline{4}}\, 1_{\underline{7}}^{\underline{3}}\, 1_{\underline{5}}^{\underline{3}}\, 1_{\underline{5}}^{\underline{2}}\, 1_{\underline{3}}^{\underline{1}}\, 1_{\underline{1}}^{\underline{1}}\, 1_{\underline{3}}^{0}\, 2_{\underline{1}}^{0}\, 1_{\underline{1}}^{1}\, 2_{1}^{1}\, 1_{3}^{2}\, 2_{3}^{3}\, 1_{5}^{3}\, 1_{5}^{4}\, 1_{7}^{4}\, 1_{9}^{6}$
}

\expandafter\def\csname KhovanovT_K15n_51748\endcsname{%
$1_{\underline{13}}^{\underline{7}}\, 1_{\underline{9}}^{\underline{5}}\, 1_{\underline{9}}^{\underline{4}}\, 1_{\underline{7}}^{\underline{4}}\, 1_{\underline{7}}^{\underline{3}}\, 1_{\underline{5}}^{\underline{2}}\, 1_{\underline{3}}^{\underline{2}}\, 1_{\underline{5}}^{\underline{1}}\, 2_{\underline{3}}^{\underline{1}}\, 1_{\underline{3}}^{0}\, 2_{\underline{1}}^{0}\, 1_{1}^{1}\, 2_{1}^{2}\, 1_{3}^{2}\, 1_{3}^{3}\, 1_{5}^{3}\, 1_{7}^{5}$
}

\expandafter\def\csname KhovanovT_K15n_57436\endcsname{%
$1_{\underline{13}}^{\underline{7}}\, 1_{\underline{11}}^{\underline{5}}\, 1_{\underline{9}}^{\underline{5}}\, 1_{\underline{9}}^{\underline{4}}\, 2_{\underline{7}}^{\underline{4}}\, 2_{\underline{7}}^{\underline{3}}\, 2_{\underline{5}}^{\underline{2}}\, 1_{\underline{3}}^{\underline{2}}\, 3_{\underline{3}}^{\underline{1}}\, 1_{\underline{1}}^{\underline{1}}\, 1_{\underline{3}}^{0}\, 2_{\underline{1}}^{0}\, 1_{\underline{1}}^{1}\, 1_{1}^{1}\, 1_{1}^{2}\, 2_{3}^{2}\, 2_{3}^{3}\, 1_{5}^{4}\, 1_{7}^{5}\, 1_{9}^{6}$
}

\expandafter\def\csname KhovanovT_K15n_57606\endcsname{%
$1_{\underline{13}}^{\underline{7}}\, 1_{\underline{9}}^{\underline{5}}\, 1_{\underline{9}}^{\underline{4}}\, 1_{\underline{7}}^{\underline{4}}\, 2_{\underline{7}}^{\underline{3}}\, 1_{\underline{5}}^{\underline{2}}\, 1_{\underline{3}}^{\underline{2}}\, 1_{\underline{5}}^{\underline{1}}\, 3_{\underline{3}}^{\underline{1}}\, 1_{\underline{3}}^{0}\, 3_{\underline{1}}^{0}\, 1_{\underline{1}}^{1}\, 1_{1}^{1}\, 2_{1}^{2}\, 2_{3}^{2}\, 2_{3}^{3}\, 1_{5}^{3}\, 1_{5}^{4}\, 1_{7}^{5}\, 1_{9}^{6}$
}

\expandafter\def\csname KhovanovT_K15n_115375\endcsname{%
$1_{\underline{13}}^{\underline{7}}\, 1_{\underline{11}}^{\underline{5}}\, 1_{\underline{9}}^{\underline{5}}\, 1_{\underline{9}}^{\underline{4}}\, 2_{\underline{7}}^{\underline{4}}\, 1_{\underline{7}}^{\underline{3}}\, 2_{\underline{5}}^{\underline{2}}\, 1_{\underline{3}}^{\underline{2}}\, 2_{\underline{3}}^{\underline{1}}\, 1_{\underline{1}}^{\underline{1}}\, 1_{\underline{3}}^{0}\, 1_{\underline{1}}^{0}\, 1_{1}^{1}\, 1_{1}^{2}\, 1_{3}^{2}\, 1_{3}^{3}\, 1_{7}^{5}$
}

\expandafter\def\csname KhovanovT_K15n_133697\endcsname{%
$1_{\underline{11}}^{\underline{6}}\, 1_{\underline{9}}^{\underline{5}}\, 1_{\underline{7}}^{\underline{4}}\, 1_{\underline{7}}^{\underline{3}}\, 2_{\underline{5}}^{\underline{3}}\, 2_{\underline{5}}^{\underline{2}}\, 2_{\underline{3}}^{\underline{2}}\, 1_{\underline{3}}^{\underline{1}}\, 1_{\underline{1}}^{\underline{1}}\, 3_{\underline{1}}^{0}\, 1_{1}^{0}\, 3_{1}^{1}\, 1_{3}^{1}\, 1_{1}^{2}\, 1_{3}^{2}\, 2_{5}^{3}\, 1_{5}^{4}\, 1_{7}^{4}\, 1_{7}^{5}\, 1_{11}^{7}$
}

\expandafter\def\csname KhovanovT_K15n_135711\endcsname{%
$1_{\underline{11}}^{\underline{6}}\, 1_{\underline{9}}^{\underline{5}}\, 1_{\underline{7}}^{\underline{4}}\, 2_{\underline{5}}^{\underline{3}}\, 2_{\underline{5}}^{\underline{2}}\, 1_{\underline{3}}^{\underline{2}}\, 1_{\underline{3}}^{\underline{1}}\, 1_{\underline{1}}^{\underline{1}}\, 2_{\underline{1}}^{0}\, 1_{1}^{0}\, 1_{\underline{1}}^{1}\, 3_{1}^{1}\, 1_{1}^{2}\, 2_{3}^{2}\, 2_{5}^{3}\, 2_{5}^{4}\, 1_{7}^{4}\, 1_{7}^{5}\, 1_{9}^{5}\, 1_{11}^{7}$
}

\expandafter\def\csname KhovanovT_K15n_148673\endcsname{%
$1_{\underline{5}}^{\underline{4}}\, 1_{\underline{3}}^{\underline{3}}\, 1_{\underline{1}}^{\underline{2}}\, 1_{\underline{1}}^{\underline{1}}\, 2_{1}^{\underline{1}}\, 2_{1}^{0}\, 2_{3}^{0}\, 1_{3}^{1}\, 1_{5}^{1}\, 3_{5}^{2}\, 1_{7}^{2}\, 4_{7}^{3}\, 1_{9}^{3}\, 1_{7}^{4}\, 1_{9}^{4}\, 2_{11}^{5}\, 1_{11}^{6}\, 1_{13}^{6}\, 1_{13}^{7}\, 1_{17}^{9}$
}

\expandafter\def\csname KhovanovT_K15n_151500M\endcsname{%
$1_{\underline{5}}^{\underline{4}}\, 1_{\underline{3}}^{\underline{3}}\, 1_{\underline{1}}^{\underline{2}}\, 2_{1}^{\underline{1}}\, 2_{1}^{0}\, 1_{3}^{0}\, 1_{3}^{1}\, 1_{5}^{1}\, 2_{5}^{2}\, 1_{7}^{2}\, 1_{5}^{3}\, 4_{7}^{3}\, 1_{7}^{4}\, 2_{9}^{4}\, 2_{11}^{5}\, 2_{11}^{6}\, 1_{13}^{6}\, 1_{13}^{7}\, 1_{15}^{7}\, 1_{17}^{9}$
}

\printdata{K14n_22185}{K14n_22589}{14^n_{22185}}{14^n_{22589}}{0}
\printdata{K15n_57436}{K15n_57606}{15^n_{57436}}{15^n_{57606}}{0}
\printdata{K15n_133697}{K15n_135711}{15^n_{133697}}{15^n_{135711}}{0}
\printdata{K15n_115375}{K15n_51748}{15^n_{115375}}{15^n_{51748 }}{0}
\printdata{K15n_148673}{K15n_151500M}{15^n_{148673}}{\overline{15}^n_{151500}}{-2}

\begin{table}[htp]
\begin{center}
\arraycolsep 2.9pt\def\arraystretch{1.2}
\rotatebox{90}{$\begin{array}{|c|rrrrrrrrrrr|}
\hline
\mathbf{K_{75}}& \mathbf{l^{-6}}& \mathbf{l^{-4}}& \mathbf{l^{-2}}& \mathbf{1}& \mathbf{l^2}&
\mathbf{l^4}& \mathbf{l^6}& \mathbf{l^{8}}& \mathbf{l^{10}}& \mathbf{l^{12}}& \mathbf{l^{14}} \\
\hline
\mathbf{1}     &       &       &     -2&     -1&      3&      3&     -1&     -1&       &       &        \\
\mathbf{m^{2}} &      7&     56&    139&    135&     25&    -23&    -11&    -22&    -21&     -5&        \\
\mathbf{m^{4}} &   -211&   -953&  -1458&   -523&    454&   -151&   -622&    -90&    128&    -14&    -26 \\
\mathbf{m^{6}} &   1579&   5441&   4719&  -3552&  -4992&   5025&   8085&   1863&   -426&    586&    365 \\
\mathbf{m^{8}} &  -5299& -14273&     77&  30645&  15926& -36529& -42026&  -8337&    844&  -4408&  -2181 \\
\mathbf{m^{10}}&   9130&  16660& -36120& -94856& -14197& 132597& 123772&  19005&    -21&  16978&   7239 \\
\mathbf{m^{12}}&  -7427&   -370&  97671& 154882& -38199&-291772&-234023& -25052&  -4083& -39284& -14827 \\
\mathbf{m^{14}}&   -161& -22770&-125295&-136384& 136281& 425309& 301818&  17636&   9947&  59198&  19943 \\
\mathbf{m^{16}}&   6309&  29798&  83571&  43754&-205222&-430970&-273994&  -1817& -12278& -60701& -18164 \\
\mathbf{m^{18}}&  -6442& -19849& -17860&  35183& 188278& 311912& 177631&  -8519&   9103&  43287&  11276 \\
\mathbf{m^{20}}&   3412&   7996& -16095& -50831&-115309&-163389& -82376&   8513&  -4248& -21639&  -4732 \\
\mathbf{m^{22}}&  -1083&  -2018&  16221&  30148&  48742&  62044&  27027&  -4244&   1254&   7544&   1313 \\
\mathbf{m^{24}}&    207&    312&  -7152& -10708& -14310& -16893&  -6106&   1265&   -227&  -1794&   -230 \\
\mathbf{m^{26}}&    -22&    -27&   1859&   2424&   2871&   3209&    901&   -228&     23&    277&     23 \\
\mathbf{m^{28}}&      1&      1&   -293&   -344&   -376&   -403&    -78&     23&     -1&    -25&     -1 \\
\mathbf{m^{30}}&       &       &     26&     28&     29&     30&      3&     -1&       &      1&        \\
\mathbf{m^{32}}&       &       &     -1&     -1&     -1&     -1&       &       &       &       &        \\
\hline
\end{array}$}
\vspace{0.5cm}
\caption{Coefficients of the HOMFLY-PT polynomials of the knot $K_{75}$}
\label{tbl:HOMFLY-PT1}
\end{center}
\end{table}

\begin{table}[htp]
\begin{center}
\arraycolsep 2.9pt\def\arraystretch{1.2}
\rotatebox{90}{$\begin{array}{|c|rrrrrrrrrrr|}
\hline
\mathbf{K_{75}^\tau}& \mathbf{l^{-6}}& \mathbf{l^{-4}}& \mathbf{l^{-2}}& \mathbf{1}& \mathbf{l^2}&
\mathbf{l^4}& \mathbf{l^6}& \mathbf{l^{8}}& \mathbf{l^{10}}& \mathbf{l^{12}}& \mathbf{l^{14}} \\
\hline
\mathbf{1}     &       &       &     -2&     -1&      3&      3&     -1&     -1&       &       &        \\
\mathbf{m^{2}} &      7&     56&    139&    135&     25&    -23&    -11&    -22&    -21&     -5&        \\
\mathbf{m^{4}} &   -211&   -964&  -1533&   -740&    111&   -466&   -783&   -125&    131&    -12&    -26 \\
\mathbf{m^{6}} &   1579&   5507&   5179&  -2207&  -2871&   6936&   9016&   2038&   -451&    577&    366 \\
\mathbf{m^{8}} &  -5299& -14405&  -1058&  27131&  10403& -41352& -44238&  -8694&    894&  -4402&  -2181 \\
\mathbf{m^{10}}&   9130&  16781& -34668& -89806&  -6231& 139236& 126581&  19388&    -56&  16977&   7239 \\
\mathbf{m^{12}}&  -7427&   -425&  96604& 150503& -45229&-297275&-236105& -25284&  -4073& -39284& -14827 \\
\mathbf{m^{14}}&   -161& -22758&-124827&-134003& 140223& 428170& 302742&  17715&   9946&  59198&  19943 \\
\mathbf{m^{16}}&   6309&  29797&  83450&  42938&-206629&-431908&-274235&  -1831& -12278& -60701& -18164 \\
\mathbf{m^{18}}&  -6442& -19849& -17843&  35354& 188587& 312100& 177665&  -8518&   9103&  43287&  11276 \\
\mathbf{m^{20}}&   3412&   7996& -16096& -50851&-115347&-163410& -82378&   8513&  -4248& -21639&  -4732 \\
\mathbf{m^{22}}&  -1083&  -2018&  16221&  30149&  48744&  62045&  27027&  -4244&   1254&   7544&   1313 \\
\mathbf{m^{24}}&    207&    312&  -7152& -10708& -14310& -16893&  -6106&   1265&   -227&  -1794&   -230 \\
\mathbf{m^{26}}&    -22&    -27&   1859&   2424&   2871&   3209&    901&   -228&     23&    277&     23 \\
\mathbf{m^{28}}&      1&      1&   -293&   -344&   -376&   -403&    -78&     23&     -1&    -25&     -1 \\
\mathbf{m^{30}}&       &       &     26&     28&     29&     30&      3&     -1&       &      1&        \\
\mathbf{m^{32}}&       &       &     -1&     -1&     -1&     -1&       &       &       &       &        \\
\hline
\end{array}$}
\vspace{0.5cm}
\caption{Coefficients of the HOMFLY-PT polynomials of the cabled mutant knot $K_{75}^\tau$.}
\label{tbl:HOMFLY-PT2}
\end{center}
\end{table}

\clearpage

\ifx\undefined\bysame
\newcommand{\bysame}{\leavevmode\hbox to3em{\hrulefill}\,}\fi

\pagebreak

\swapnumbers
\theoremstyle{plain}
\newtheorem*{thmjones}{Theorem 3.2}
\newtheorem*{thmalex}{Theorem 2.9}

\begin{center}
\Large \textbf{Erratum added February 2012}
\end{center}

\begin{abstract}
  Proposition 2.7 of the original paper \cite{DGST} is false and as a
  result Corollary~2.8 has not been established.  Here, we provide
  alternate proofs of the results in our paper which depended on those
  claims, with the exception of the invariance of generalized knot
  signatures.  In particular, all the results claimed in Table 1.2 of
  the original paper have still been proved.
\end{abstract}

\maketitle

The two places where Corollary 2.8 was used are Theorems 2.9 and
3.2.  We start by giving a correct proof of Theorem 3.2. 

\begin{thmjones}
The colored Jones polynomials of a knot are invariant under $(2,0)$-mutation
for all colors.  
\end{thmjones}

\begin{proof}[Proof of Theorem 3.2]
  Let $F$ be a closed genus 2 surface in $S^3$ disjoint from a knot
  $K$, and let $K^\tau$ be the mutant of $K$ along $F$, where here
  $\tau$ is the hyperelliptic involution.  We will use that
  the colored Jones polynomials can be defined via the Kauffman
  bracket skein module (KBSM), in the style of topological quantum
  field theory.

  The key here is that by Theorem 3.1 of \cite{P}, one has the
  following basis for the KBSM of $F \times I$ where $I = [-1,1]$: the
  set of isotopy classes of unoriented links in $F \times \{0\}$ where
  every component of the link is an essential curve.  Here, each such
  curve is given the blackboard framing.  Now the hyperelliptic
  involution $\tau$ acts trivially on this set of framed links and
  therefore also on $\mathrm{KBSM}(F \times I)$.
  
  The surface $F$ divides $S^3 \setminus K$ into two pieces, which we
  denote by $X$ and $Y$.  Then $S^3\setminus K^{\tau}$ is obtained by
  gluing $X$ to one side of $F \times I$ and $Y$ to the other side via
  the hyperelliptic involution $\tau$.  As $\tau$ acts trivially on
  $\mathrm{KBSM}(F \times I)$, it follows that
  $\mathrm{KBSM}(S^3\setminus K)$ is isomorphic to
  $\mathrm{KBSM}(S^3\setminus K^{\tau})$.  By Masbaum and Vogel
  \cite{MV}, it follows that the colored Jones polynomials of $K$ and
  $K^{\tau}$ are equal for all colors.
\end{proof}

We next give a correct proof of part of Theorem 2.9. 
 
\begin{thmalex}[Revised]
  The Alexander polynomial of a knot in $S^3$ does not change under
  $(2,0)$-mutation.
\end{thmalex}

The statement of Theorem 2.9 in \cite{DGST} asserts that the
generalized signatures are also invariant under $(2,0)$-mutation, but
we do not know how to establish this; these signatures are invariant
under genus 2 \emph{handlebody} mutation, see \cite{CL}.

\begin{proof}  
  The Alexander polynomial of a knot is determined by all of its
  colored Jones polynomials (this is the Melvin-Morton-Rozansky
  Conjecture, which was proven in~\cite{B-NG}).  Thus Theorem 3.2
  implies that the Alexander polynomials does not change under
  $(2,0)$-mutation.
\end{proof}

\subsection*{The problem with Proposition 2.7} 

Proposition 2.7 claimed that if $K$ is a knot in $S^3$ which is
disjoint from a genus 2 surface $F$, then either $K^\tau$ is obtained
from $K$ by various kinds of \emph{handlebody} mutation or $K^\tau
\cong K$.  In particular, we claimed that if $F$ is incompressible in
the complement of $K$, then in fact $F$ bounds a handlebody in $S^3$;
this is simply false, as the following example shows.   Start with a
knotted solid torus $V$ in $S^3$.   If we then drill out a tunnel
from $V$, we get a submanifold $Y$ with $F = \partial X$ a genus 2
surface; by choosing a complicated tunnel, we can arrange that
$F$ is incompressible in $Y$.   Let $X$ be the complement of $Y$, and
choose a knot $K$ in $X$ which runs through the tunnel and is chosen
so that $F$ is incompressible in $X \setminus K$.   Then $F$ is
incompressible in $S^3 \setminus K$, but it does not bound a
handlebody on either side; hence mutation along $F$ is \emph{not}
$(2,0)$--handlebody mutation.  
 
\subsection*{Acknowledgment}

We are extremely grateful to Mario Eudave Mu\~noz for finding the error in
Proposition 2.7 and providing the above
counter\hyp example.

\end{document}

%% file: pictures/symsurfaces.pstex_t
\begin{picture}(0,0)%
\includegraphics{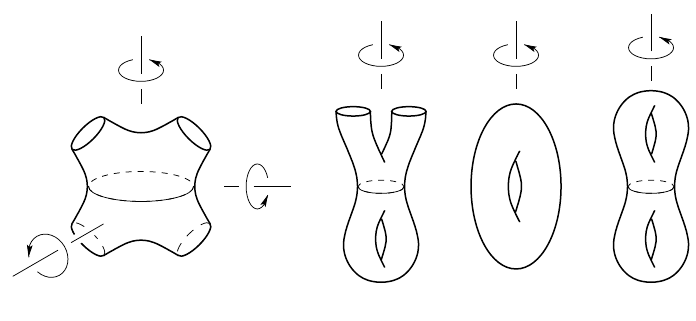}%
\end{picture}%
\setlength{\unitlength}{3158sp}%
\begingroup\makeatletter\ifx\SetFigFont\undefined%
\gdef\SetFigFont#1#2#3#4#5{%
  \reset@font\fontsize{#1}{#2pt}%
  \fontfamily{#3}\fontseries{#4}\fontshape{#5}%
  \selectfont}%
\fi\endgroup%
\begin{picture}(6956,3197)(165,-2680)
\put(2851,-1711){\makebox(0,0)[lb]{\smash{{\SetFigFont{10}{12.0}{\rmdefault}{\mddefault}{\updefault}{\color[rgb]{0,0,0}$\tau_2$}%
}}}}
\put(826,-2386){\makebox(0,0)[lb]{\smash{{\SetFigFont{10}{12.0}{\rmdefault}{\mddefault}{\updefault}{\color[rgb]{0,0,0}$\tau_1$}%
}}}}
\put(1876,-211){\makebox(0,0)[lb]{\smash{{\SetFigFont{10}{12.0}{\rmdefault}{\mddefault}{\updefault}{\color[rgb]{0,0,0}$\tau_3$}%
}}}}
\put(4276,-61){\makebox(0,0)[lb]{\smash{{\SetFigFont{10}{12.0}{\rmdefault}{\mddefault}{\updefault}{\color[rgb]{0,0,0}$\tau$}%
}}}}
\put(5626,-61){\makebox(0,0)[lb]{\smash{{\SetFigFont{10}{12.0}{\rmdefault}{\mddefault}{\updefault}{\color[rgb]{0,0,0}$\tau$}%
}}}}
\put(6976, 14){\makebox(0,0)[lb]{\smash{{\SetFigFont{10}{12.0}{\rmdefault}{\mddefault}{\updefault}{\color[rgb]{0,0,0}$\tau$}%
}}}}
\put(3976,-2611){\makebox(0,0)[b]{\smash{{\SetFigFont{10}{12.0}{\rmdefault}{\mddefault}{\updefault}{\color[rgb]{0,0,0}Type (1,2)}%
}}}}
\put(5326,-2611){\makebox(0,0)[b]{\smash{{\SetFigFont{10}{12.0}{\rmdefault}{\mddefault}{\updefault}{\color[rgb]{0,0,0}Type (1,0)}%
}}}}
\put(6676,-2611){\makebox(0,0)[b]{\smash{{\SetFigFont{10}{12.0}{\rmdefault}{\mddefault}{\updefault}{\color[rgb]{0,0,0}Type (2,0)}%
}}}}
\put(1576,-2611){\makebox(0,0)[b]{\smash{{\SetFigFont{10}{12.0}{\rmdefault}{\mddefault}{\updefault}{\color[rgb]{0,0,0}Type (0,4)}%
}}}}
\end{picture}%

%% file: pictures/skein-generators.pstex_t
\begin{picture}(0,0)%
\includegraphics{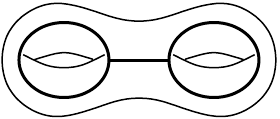}%
\end{picture}%
\setlength{\unitlength}{3158sp}%
\begingroup\makeatletter\ifx\SetFigFont\undefined%
\gdef\SetFigFont#1#2#3#4#5{%
  \reset@font\fontsize{#1}{#2pt}%
  \fontfamily{#3}\fontseries{#4}\fontshape{#5}%
  \selectfont}%
\fi\endgroup%
\begin{picture}(2780,1176)(2211,-1548)
\put(3601,-886){\makebox(0,0)[b]{\smash{{\SetFigFont{10}{12.0}{\rmdefault}{\mddefault}{\updefault}{\color[rgb]{0,0,0}b}%
}}}}
\put(4351,-1261){\makebox(0,0)[b]{\smash{{\SetFigFont{10}{12.0}{\rmdefault}{\mddefault}{\updefault}{\color[rgb]{0,0,0}c}%
}}}}
\put(2851,-1261){\makebox(0,0)[b]{\smash{{\SetFigFont{10}{12.0}{\rmdefault}{\mddefault}{\updefault}{\color[rgb]{0,0,0}a}%
}}}}
\end{picture}%

%% file: pictures/tangles-cabled.pstex_t
\begin{picture}(0,0)%
\includegraphics{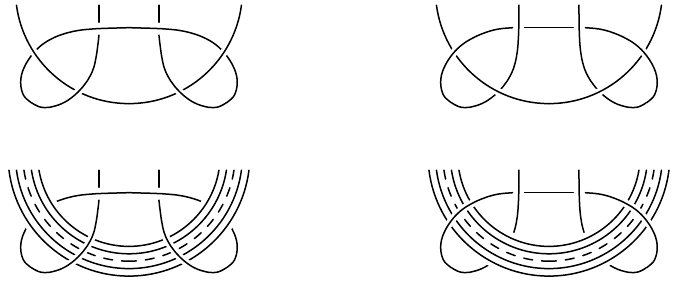}%
\end{picture}%
\setlength{\unitlength}{3158sp}%
\begingroup\makeatletter\ifx\SetFigFont\undefined%
\gdef\SetFigFont#1#2#3#4#5{%
  \reset@font\fontsize{#1}{#2pt}%
  \fontfamily{#3}\fontseries{#4}\fontshape{#5}%
  \selectfont}%
\fi\endgroup%
\begin{picture}(6750,2806)(1036,-3122)
\put(1051,-961){\makebox(0,0)[rb]{\smash{{\SetFigFont{10}{12.0}{\rmdefault}{\mddefault}{\updefault}{\color[rgb]{0,0,0}$T$:}%
}}}}
\put(5251,-961){\makebox(0,0)[rb]{\smash{{\SetFigFont{10}{12.0}{\rmdefault}{\mddefault}{\updefault}{\color[rgb]{0,0,0}$T^\tau$:}%
}}}}
\put(7576,-2161){\makebox(0,0)[b]{\smash{{\SetFigFont{10}{12.0}{\rmdefault}{\mddefault}{\updefault}{\color[rgb]{0,0,0}\vbox{\halign{\hfill#\hfill\cr$n$ strands\cr$\overbrace{\strut\hskip 0.3truein}$\cr}}}%
}}}}
\put(5251,-2611){\makebox(0,0)[rb]{\smash{{\SetFigFont{10}{12.0}{\rmdefault}{\mddefault}{\updefault}{\color[rgb]{0,0,0}$T^\tau(1,n)$:}%
}}}}
\put(1051,-2611){\makebox(0,0)[rb]{\smash{{\SetFigFont{10}{12.0}{\rmdefault}{\mddefault}{\updefault}{\color[rgb]{0,0,0}$T(1,n)$:}%
}}}}
\put(3376,-2161){\makebox(0,0)[b]{\smash{{\SetFigFont{10}{12.0}{\rmdefault}{\mddefault}{\updefault}{\color[rgb]{0,0,0}\vbox{\halign{\hfill#\hfill\cr$n$ strands\cr$\overbrace{\strut\hskip 0.3truein}$\cr}}}%
}}}}
\end{picture}%

%% file: pictures/cabled-mutants.pstex_t
\begin{picture}(0,0)%
\includegraphics{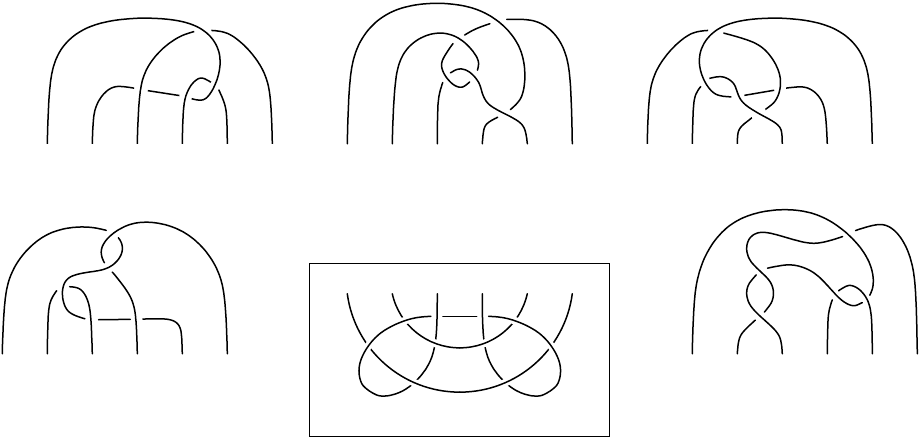}%
\end{picture}%
\setlength{\unitlength}{3158sp}%
\begingroup\makeatletter\ifx\SetFigFont\undefined%
\gdef\SetFigFont#1#2#3#4#5{%
  \reset@font\fontsize{#1}{#2pt}%
  \fontfamily{#3}\fontseries{#4}\fontshape{#5}%
  \selectfont}%
\fi\endgroup%
\begin{picture}(9194,4363)(129,-4498)
\put(1726,-1861){\makebox(0,0)[b]{\smash{{\SetFigFont{12}{12.0}{\rmdefault}{\mddefault}{\updefault}{\color[rgb]{0,0,0}a. $14^n_{22185}$ and  $14^n_{22589}$.}%
}}}}
\put(4726,-4336){\makebox(0,0)[b]{\smash{{\SetFigFont{12}{12.0}{\rmdefault}{\mddefault}{\updefault}{\color[rgb]{0,0,0}Common closure: $T(1,2)$}%
}}}}
\put(7726,-1861){\makebox(0,0)[b]{\smash{{\SetFigFont{12}{12.0}{\rmdefault}{\mddefault}{\updefault}{\color[rgb]{0,0,0}c. $15^n_{115375}$ and  $15^n_{51748}$.}%
}}}}
\put(4726,-1861){\makebox(0,0)[b]{\smash{{\SetFigFont{12}{12.0}{\rmdefault}{\mddefault}{\updefault}{\color[rgb]{0,0,0}b. $15^n_{57606}$ and  $15^n_{57436}$.}%
}}}}
\put(8176,-3961){\makebox(0,0)[b]{\smash{{\SetFigFont{12}{12.0}{\rmdefault}{\mddefault}{\updefault}{\color[rgb]{0,0,0}e. $15^n_{148673}$ and  $\overline{15}^n_{151500}$.}%
}}}}
\put(1276,-3961){\makebox(0,0)[b]{\smash{{\SetFigFont{12}{12.0}{\rmdefault}{\mddefault}{\updefault}{\color[rgb]{0,0,0}d. $15^n_{133697}$ and  $15^n_{135711}$.}%
}}}}
\end{picture}%